\documentclass[11pt]{article}

\usepackage{latexsym}
\usepackage{amsmath,amssymb}
\usepackage{graphicx}
\usepackage{psfrag}
\newtheorem{Theorem}{Theorem}[part]
\newtheorem{Definition}{Definition}[part]

\newtheorem{Assumption}{Assumption}[part]
\newtheorem{Lemma}{Lemma}[part]
\newtheorem{Corollary}{Corollary}[part]
\newtheorem{Remark}{Remark}[part]


\def \Sum{\displaystyle\sum}

\def \Frac{\displaystyle\frac}
\def \Inf{\displaystyle\inf}
\def \Sup{\displaystyle\sup}
\def \Lim{\displaystyle\lim}
\def \Liminf{\displaystyle\liminf}
\def \Limsup{\displaystyle\limsup}
\def \Max{\displaystyle\max}
\def \Min{\displaystyle\min}

\def \N{I\!\!N}
\def \R{I\!\!R}

\def \F{I\!\!F}

\def \Ac{{\cal A}}

\def \Cc{{\cal C}}
\def \Dc{{\cal D}}
\def \Ec{{\cal E}}
\def \Fc{{\cal F}}

\def \Ic{{\cal I}}
\def \Lc{{\cal L}}
\def \Pc{{\cal P}}

\def \Sc{{\cal S}}
\def \Tc{{\cal T}}

\def \Xc{{\cal X}}
\def \Yc{{\cal Y}}

\def \ep{\hbox{ }\hfill$\Box$}

\def\Dt#1{\Frac{\partial #1}{\partial t}}

\def\Dy#1{\Frac{\partial #1}{\partial y}}

\def\reff#1{{\rm(\ref{#1})}}

\begin{document}
\title{Optimal insurance demand under marked point processes shocks: a dynamic
  programming duality approach}
\author{Mohamed MNIF  \\
  \small LAMSIN\\
            \small Ecole Nationale d'Ing\'enieurs de Tunis \\
            \small B.P. 37,
            1002, Tunis Belv\'ed\`ere, Tunisie\\
 mohamed.mnif@enit.rnu.tn }
\date{August 18, 2010}
\maketitle

\begin{abstract}
We study the stochastic control problem of maximizing expected utility from
terminal wealth under a non-bankruptcy constraint. The wealth process is subject to shocks produced by a 
general marked point process. The problem of the agent is to derive the
optimal insurance strategy which allows "lowering" the level of the shocks.
This optimization problem is related to a
suitable dual stochastic control problem in which the delicate boundary constraints disappear. 
We characterize the dual value function as the unique viscosity solution of the corresponding
a Hamilton Jacobi Bellman Variational Inequality (HJBVI in short). 
\end{abstract}

\vspace{7mm}

\noindent {\bf Key words~:} Optimal insurance; stochastic control; duality; optional decomposition; dynamic  programming principle; viscosity solution

\vspace{5mm}

\noindent {\bf MSC Classification (2000)~:} 
93E20, 60J75, 49L25.
\newpage
\section{Introduction}

\setcounter{equation}{0}
\setcounter{Assumption}{0}
\setcounter{Example}{0}
\setcounter{Theorem}{0}
\setcounter{Proposition}{0}
\setcounter{Corollary}{0}
\setcounter{Lemma}{0}
\setcounter{Definition}{0}
\setcounter{Remark}{0}
We study the optimal insurance demand problem of an agent whose wealth is
subject to shocks produced by some marked point process. Such a problem was
formulated by Bryis \cite{bri86} in continuous-time with Poisson shocks. 
Gollier \cite{gol94} studied a similar problem where shocks are not proportional to
wealth. An explicit solution to the problem is provided by Bryis by writing
the associated Hamilton-Jacobi-Bellman (HJB in short) equation. In Bryis \cite{bri86} and
Gollier \cite{gol94}, they modeled the insurance premium by an affine function of the
insurance strategy $\theta=(\theta_t)_{t\in [0,T]}$ which is the
rate of insurance decided to be covered by the agent. If the agent is
subject to some accident at time $t$ which costs an amount $Z$, then he will pay $\theta_t Z$  
and the insurance company reimburses the amount
$(1-\theta_t)Z$. They didn't assume any constraint on the
insurance strategy which is not realistic. \\
In risk theory, Hipp and Plum \cite {hipplu00} analysed the trading
strategy, in risky  assets, which is optimal with respect to the criterion of
minimizing the ruin 
probability. They derived the HJB equation related
to this problem and proved the existence of a solution and a verification
theorem. When the claims are exponentially distributed, the ruin probability
decreases exponentially and the optimal amount invested in risky assets
converges to a constant independent of the reserve level. Hipp and Schmidli
\cite {hipplu01} have obtained the asymptotic behaviour of the ruin
probability under the optimal investment strategy in the small claim
case. Schmidli \cite{sch99}  studied the optimal proportional
reinsurance policy which minimizes the ruin probability in infinite horizon.
He derived the associated HJB equation, proved the
existence of a solution and a verification theorem in the diffusion case. He
proved that the ruin probability decreases exponentially whereas the optimal
proportion to insure is constant. Moreover, he gave some conjecture in the
Cram\'er-Lundberg case. H\o jgaard and Taksar \cite {hoj98} studied another problem
of proportional reinsurance. They considered the issue of reinsurance optimal
fraction, that maximizes the return function. They modelled the reserve  process as  a
diffusion process. \\ 
Touzi \cite{tou00} studied the problem of maximizing the expected utility from
terminal wealth when the insurance strategy 
is valued in $[0,1]$ at each time . He modeled the wealth process by a
Dol\'eans-Dade exponential process. He 
assumed a  boundedness assumption on  
the jump term which guarantees the positivity of the wealth process. He solved
this stochastic control problem by using duality method. \\ 
Duality method was introduced by Karatzas et al. \cite{karle} and Cox and
Huang \cite{cox89}. Cox and Huang characterized the optimal consumption- portfolio
policies when there exist non-negativity constraints on consumption and on final
wealth. They gave a verification theorem which involves a linear partial
differential equation unlike the nonlinear Bellman equation. In few cases they constructed
the optimal control. Extensions to the case of constrained investment are
considered by Cvitani\'c and Karatzas \cite{cvikar92} and to the case of incomplete
markets by Karatzas et al. \cite{karlehshrxu91}. Typically, in
incomplete markets, we have to solve a dual problem which leads in the Markov
case to nonlinear partial differential equation.\\
In this paper, we  model the claims by using a compound Poisson process. 
The insurance trading strategy is constrained to remain in $[0,1]$. We
impose a constraint of non-bankruptcy on the wealth process $X_t$ of the
agent for all $t$. 
The objective of the agent is to maximize the expected utility of
the terminal wealth over all admissible strategies and to determine  the
optimal policy of insurance. \\
In our case the wealth process positivity constraint is a real one unlike 
the problem formulated in Touzi \cite{tou00}. \\
Our stochastic optimization problem is a particular
case of a general structure of problems developed in Mnif and Pham \cite{mnipham01} 
who considered the following optimization problem:
\begin{eqnarray} \label{optimdyn}
\Max_{X\in\Xc_+(x)} E[U(X_T)], \;\;\; x\in \R, 
\end{eqnarray}  
where  $\Xc_+(x):=\left\{ x+ X~: X \in \Xc \mbox{ s.t }\,\, X_t\geq 0 \mbox{
    for all } 0\leq t\leq T\right\}$ where $\Xc$ is a family of semi-martingales. 
Existence and uniqueness of
solution of problem (\ref{optimdyn}) is then proved. The
optimal solution characterization is obtained from a dual formulation under
    minimal assumptions on the objective function.\\
In this paper, we study the dual value function by a PDE approach. The dual problem
    appears as a mixed control/singular optimization problem with dynamics
    $(Y_t=Z_t D_t,\,t\in [0,T])$ governed by a classical control term $Z$
    which comes from the insurance strategy and a singular term $D$ which comes 
from the
    state constraint.\\   
 The originality of this paper is to study a stochastic control problem with state constraint. 
The wealth of the investor must be non-negative even after a jump. 
The duality method is not another alternative to solve this problem. In fact the
    primal problem leads to a HJB equation with boundary conditions. 
    Because of the  
state space constraints these boundary conditions are not obvious to obtain. 
However, these delicate boundary conditions disappear in the dual problem.
The regularity of the dual value function is not obvious to obtain. 
This explains the use of the notion of discontinuous viscosity solutions. 
The comparison theorem is not proved in a general framework since the operator 
which appears in the HJBVI contains an inf on a unbounded set which makes it 
discontinuous. In this paper, we prove the comparison theorem 
only when the space of claims is a finite one.\\
The paper is organized as follows. Section 2 describes the
model. In Section 3, we  formulate the dual optimization problem and we derive
the associated HJBVI 
for the value function. In Section 4, we prove that the
    dual value function is 
a viscosity solution of our HJBVI. In Section 5, we
prove a comparison theorem.

\section{Problem formulation}

\setcounter{equation}{0}
\setcounter{Assumption}{0}
\setcounter{Example}{0}
\setcounter{Theorem}{0}
\setcounter{Proposition}{0}
\setcounter{Corollary}{0}
\setcounter{Lemma}{0}
\setcounter{Definition}{0}
\setcounter{Remark}{0}

Let $(\Omega,\Fc,P)$ be a complete probability space. 
We assume that the claims are generated by a compound Poisson process. 
More precisely, we consider an integer-valued random measure $\mu (dt,dz)$
with  compensator $\pi(dz)dt$. We  assume that $\pi(dz)=\varrho G(dz)$ where
$G(dz)$ is a probability distribution on the bounded set $C\subseteq \R_+$ and
$\varrho$ is a 
positive constant. In this case, the integral, with respect to the random  
measure $\mu (dt,dz)$, is  simply a  compound Poisson process: we have $ 
\int_0^t \int_C z\mu(du,dz)=\sum_{i= 1}^{N_t}Z_i$, where $N=\{N_t,t\geq 0\}$ 
is a Poisson process with intensity $\varrho$ and $\{Z_i,i\in \N\}$ is a 
sequence of random variables with common distribution $G$ which represent the
claim sizes.\\   
Let $T>0$ be a finite time horizon. We denote by $\F=(\Fc_t)_{0\leq t\leq T}$ 
the filtration generated by the random  
measure $\mu (dt,dz)$.\\ 
By definition of the intensity $\pi(dz)dt$, the compensated jump process:
\begin{eqnarray*}
\tilde \mu(dt,dz):= \mu(dt,dz)-\pi(dz)dt
\end{eqnarray*} 
is such that $\{\tilde \mu([0,t]\times B),0\leq t\leq T\}$ is a $(P,\F)$
martingale for all $B\in \Cc$, where $\Cc$ is the Borel
$\sigma$-field on $C$. \\
An insurance strategy is a predictable process $\theta=(\theta_t)_{0\leq t\leq
  T}$ which represents the rate of insurance covered by the agent.
We assume that the insurance premium is an affine function of the insurance
  strategy. 
Given an initial wealth $x\geq 0$ at time t and an insurance strategy
$\theta$, the wealth process of the agent at time $s\in [t,T]$ is then given by~:
\begin{eqnarray}\label{eds}
X_s^{t,x,\theta} &:=& x + \int_t^s \left( \alpha - \beta (1-\theta_u)
\right) du - \int_t^s \int_C\theta_u  z \mu(du,dz).
\end{eqnarray}
We assume that $\alpha\geq \beta \geq 0$ which means that the premium rate received
by the agent is lower then the premium rate paid to the insurer. 
In the literature, this problem is known as a proportional reinsurance one. 
The agent is an insurer who has to pay a premium to the reinsurer.
We impose that the insurance strategy satisfies:
\begin{eqnarray}\label{bornekm}
\theta_s \in [0,1]\,\,\, \mbox{ a.s. for all } t\leq s \leq T.
 \end{eqnarray}
We also impose the following non-bankruptcy constraint on the wealth process:
\begin{eqnarray}\label{contraintem}
X_s^{t,x,\theta} \geq  0 \,\,\mbox{ a.s. for all } t \leq s\leq T.
\end{eqnarray}    
Given an initial wealth $x\geq 0$ at time $t$, an admissible policy $\theta$ is a
predictable stochastic process
$(\theta_s)_{t\leq s\leq T}$ , such that conditions (\ref{bornekm}) and  (\ref{contraintem}) 
are  satisfied.  We denote by
$\Ac(t,x)$ the set of all admissible policies and
$\Sc(t,x):=\{X^{t,x,\theta}\mbox{ such that } \theta\in \Ac(t,x)\}$. \\ 
Our agent has preferences modeled by a utility function $U$. \\
We assume that the agent's utility is described by a CRRA utility 
function i.e. $U(x)=\frac{x^\eta}{\eta}$, where $\eta\in (0,1)$.\\
We denote by $I$ the inverse of $U^{\prime}$ and we
introduce the conjugate function of $U$ defined by
\begin{eqnarray} \label{deftildeU}
\tilde U(y)&:=&\sup_{x>0}\{U(x)-xy\},\;\;\; y >0 \nonumber\\
&=&U(I(y))-yI(y).
\end{eqnarray}
A straightforward calculus shows that $\displaystyle{\tilde U(y)=\frac{y^{-\gamma}}{\gamma}}$ 
where $\gamma=\frac{\eta}{1-\eta}$ 
and $\tilde U'(y) = -I(y)$ for all $ y >0$.\\
The objective of the agent is to find the value function which is defined as
\begin{eqnarray}\label{pvaleur} 
v(t,x):=\Sup_{\theta \in \Ac(t,x)}E(U(X_T^{t,x,\theta})). 
\end{eqnarray}
{\bf Notations:} The constants which appear in the paper are generic and could change from line to line. 
\section{Dual optimization problem}

\setcounter{equation}{0}
\setcounter{Assumption}{0}
\setcounter{Example}{0}
\setcounter{Theorem}{0}
\setcounter{Proposition}{0}
\setcounter{Corollary}{0}
\setcounter{Lemma}{0}
\setcounter{Definition}{0}
\setcounter{Remark}{0}

First we introduce some notations. Let $x\geq 0$ and $t\in[0,T]$. We denote by
$\Pc(\Sc(t,x))$ the set of all probability measures $Q$ $\sim$ $P$ with the
following property: there exists $A$ $\in$ $\Ic_p$, set of non-decreasing
predictable processes with $A_0$ $=$ $0$, such that~: 
\begin{eqnarray} \label{pc}  
X - A  \; \mbox{is a} \; Q-\mbox{local super-martingale for any} \; 
X \in \Sc(t,x).  
\end{eqnarray} 
The upper variation process of $\Sc(t,x)$ under $Q$ $\in$ $\Pc(\Sc(t,x))$ 
is the element
$\tilde{A}^{\Sc(t,x)}(Q)$ in $\Ic_p$ satisfying \reff{pc} and such that 
$A-{\tilde A^{\Sc(t,x)}}(Q)$ $\in$
$\Ic_p$ for any $A$ $\in$ $\Ic_p$ satisfying \reff{pc}. \\
From Lemma 2.1 of  F\"ollmer and Kramkov \cite {folkab97}, we can derive
$\Pc(\Sc(t,x))$ and ${\tilde A^{\Sc(t,x)}}(Q)$. This result states that
$Q\in\Pc(\Sc(t,x))$ iff there 
is an upper bound for all the  predictable processes arising in the
Doob-Meyer decomposition of the special semi-martingale $V$ $\in$
$\Sc(t,x)$ under $Q$. In this case, the upper variation process is 
equal to this upper bound. \\
It is well-known from the martingale representation
theorem for random measures (see e.g. Br\'emaud \cite{bre81}) that all
probability measures $Q$ $\sim$ $P$ have a density process in the form~:  
\begin{eqnarray}\label{exporetiel} 
Z^\rho_s &=& \Ec \left( \int_t^s\int_C (\rho_u(z) -1) \tilde \mu(du,dz) \right),\,\, s\in [t,T],
\end{eqnarray}
where 
$\rho$ $\in$ ${\cal U}_t$ $=$ $\{(\rho_s(z))_{t\leq s\leq T}$ predictable process~:  
$\rho_s(z)$ $>$ $0$, a.s., $t\leq s \leq T$,$z\in C$,
$\int_t^T \int_C\Big(|\log\rho_s(z)| + \rho_s(z) \pi(dz)\Big) ds <\infty$  and $E[Z_T^\rho] = 1\}$.\\
By Girsanov's theorem, the predictable compensator of an element
$X^\theta \in \Sc(t,x)$ under $P^\rho$ $=$ $Z_T^\rho.P$ is~:
\begin{eqnarray*}
A_s^{\rho,\theta} &=&  \int_t^s (\alpha -\beta)du + \int_t^s \theta_u (\beta-\int_C\rho_u(z)\, z\,\pi(dz)) du.
\end{eqnarray*}  
We deduce from   Lemma 2.1 of  F\"ollmer and Kramkov \cite{folkab97} that 
$\Pc(\Sc(t,x))$ $=$ $\{P^\rho~: \rho \in {\cal U}_t\}$ and the upper
variation process of $P^\rho$ is~:
\begin{eqnarray*}
\tilde A^{\Sc(t,x)}_s(P^\rho) &=& \int_t^s (\alpha -\beta)du +\int_t^s  (\beta-\int_C\rho_u(z)\,z\,\pi(dz))_+ du.
\end{eqnarray*}       
From the non-decreasing property of $U$, we have
\begin{eqnarray*}
v(t,x)=\Sup_{H\in\Cc_+(t,x)}E[U(H)],
\end{eqnarray*}   
where $\Cc_+(t,x)=\{H\in L^0_+(\Fc_T): X_T^{t,x,\theta} \geq H \,a.s. \mbox{ for
  } 
\theta\in \Ac(t,x)\}$. It is easy to check the convexity property of the
family $\Cc_+(t,x)$. For the closure property of this family in the 
semi-martingale topology, we refer to Pham \cite{pham00}.\\  
The semi-martingale topology is associated to the Emery distance
between two semi-martingales $\tilde X^1$ and $\tilde X^2$ defined as~: 
\begin{eqnarray*} \label{emery}
D_E(\tilde X^1,\tilde X^2) &:=& \sum_{n\geq 1} 2^{-n} 
E\left[\sup_{0\leq t\leq T\wedge n}|\tilde X_t^1-\tilde X_t^2|\wedge 1 \right].
\end{eqnarray*}
We refer to M\'emin \cite{m80} for details on the semi-martingale
topology. Since $\Cc_+(t,x)$ is convex and closed and  
using the optional decomposition under constraints of
F\"ollmer and Kramkov \cite{folkab97}, Mnif and Pham \cite{mnipham01} gave the following 
dual characterization of the set  $\Cc_+(t,x)$
\begin{eqnarray}\label{carac}
& &H\in \Cc_+(t,x)\\ 
&\Longleftrightarrow& 
J(H):=\Sup_{ Z\in \Pc^0(t,x)\,,\tau\in \Tc_t}  
E\left[Z_TH  1_{\tau = T} - \int_t^\tau Z_u (\alpha -\beta +(\beta -\int_C\rho_u(z)\,z\, \pi(dz))_{+})du
\right]
\leq x\nonumber,
\end{eqnarray} 
where $\Pc^0(t,x)$ is the subset of elements $P^\rho \in \Pc(\Sc(t,x))$ such
that $\tilde A_T^{\Sc(t,x)}(P^\rho)$ is bounded and $\Tc_t$ is the set of all stopping times valued in $[0,T]$.\\
As a corollary (see their corollary 4.1 ), they deduce that the set of
admissible insurance strategies $\Ac(0,x)$ is non empty iff $ x\geq
b:=(\beta-\alpha)T$. 
\begin{Remark}
If the agent initial wealth is equal to $ x=b$ and since $\Ac(0,x)$ is
not empty, then the wealth process is given by
$X_t^{0,x,\theta}=(\beta-\alpha)(T-t)$ for all  $0\leq t\leq T$ and so the only
admissible strategy is $\theta_t=0$ for all  $0\leq t\leq T$ which implies
that the dynamic version of the value function satisfies
$v\left(t,(\beta-\alpha)(T-t)\right)=0$ for all  $0\leq t\leq T$. These
boundary conditions obtained from the duality approach are not obvious from
the primal approach. 
\end{Remark}
Now, we fix some initial wealth $x\geq b$. We make the following assumption  
\begin{Assumption} \label{hyputil}
We assume that there exist $\bar \gamma \in (0,1)$,  $\bar Q \in\Pc^0(t,x)$ 
with density $\bar Z_T = \frac {d\bar Q}{dP}$ satisfying 
\begin{eqnarray*}\label{integrable}
(\bar Z_T)^{-1}\in L^{\bar p}(P)
\end{eqnarray*}
for some $\bar p > \frac{\bar\gamma}{1-\bar\gamma}$. 
\end{Assumption}
Under Assumption \ref{hyputil}, Existence and uniqueness of problem
\reff{pvaleur}  are proved (see Mnif and Pham \cite{mnipham01}). We focus now
on the study of the dual formulation.\\
The two following lemmas allow us to give an expression of the dual value function.
We denote by
${\cal D}_t$  the set of nonnegative, 
nonincreasing predictable and c\`adl\`ag processes 
$D$ $=$ $(D_s)_{t\leq s\leq T}$ with $D_t$ $=$ $1$,
\begin{eqnarray*}
\Yc^0(t):=\{Y^{\rho,D}=Z^\rho D,\,Z^\rho \in \Pc^0(t,x),\, D\in {\cal D}_t\}
\end{eqnarray*}
for all $ t\in [0,T]$ and $L_+^0(\Fc_T)$ is the set of  nonnegative  $\Fc_T$-measurable random variables.
\begin{Remark}
{\rm We omit the dependence of $\Yc^0(t)$ in the initial wealth $x$, since $x$ is fixed in all the
paper. }
\end{Remark}
\begin{Remark}
{\rm The set ${\cal D}_t$ is introduced in the paper of Elkaroui and Jeanblanc \cite{elkjea98}
in the continuous case. 
In our case, we must enlarge this set to include c\`adl\`ag processes and extend some results of Mnif and Pham \cite{mnipham01}.}
\end{Remark}
\begin{Lemma}\label{lemdc}
For all $x$ $\geq$ $b$, $X$ $\in$ $\Sc(t,x)$, $Y$ $=$ $ZD$, $Z$ $\in$ 
$\Pc^0(t,x)$, $D$ $\in$ $\Dc_t$, the processes~:
\begin{eqnarray*}
& & Z_.X_.- \int_t^. Z_u(\alpha -\beta +(\beta -\int_C\rho_u(z)\,z\, \pi(dz))_{+})du, \\
& \mbox{and} & Y_.X_.- \int_t^. Y_u(\alpha -\beta +(\beta -\int_C\rho_u(z)\,z\, \pi(dz))_{+})du
\end{eqnarray*}
are supermartingales under $P$. 
\end{Lemma}
{\bf Proof.}
By definition of $\Pc^0(t,x)$, the process 
$Z_.(X_.- \int_t^.(\alpha -\beta +(\beta -\int_C\rho_u(z)\,z\, \pi(dz))_{+})du)$ 
is a $P$-local supermartingale. 
From Theorem VII.35 in 
Dellacherie and Meyer \cite{delmey82}, the process 
\begin{eqnarray} \label{M}
M= Z_.X_.- \int_t^.Z_u(\alpha -\beta +(\beta -\int_C\rho_u(z)\,z\, \pi(dz))_{+})du
\end{eqnarray}
is a 
$P$-local supermartingale. 
Moreover, $M$ is bounded from below by the random variable $-\int_t^.Z_u(\alpha -\beta +(\beta -\int_C\rho_u(z)\,z\, \pi(dz))_{+})du)$, which
is integrable
under $P$. We deduce by Fatou's lemma that
$M$  is a  $P$-supermartingale. 
On the other hand, by It\^o's product rule and since $D$ is predictable
with finite variation, we get~:
\begin{eqnarray} \label{surmarD}
& &Y_TX_T- \int_t^T Y_u(\alpha -\beta +(\beta -\int_C\rho_u(z)\,z\, \pi(dz))_{+})du\nonumber\\
&=& Y_tX_t + \int_t^T D_{u^-} dM_u+ 
\int_t^T  Z_{u^-} X_{u^-} dD_u + \Sum_{t\leq u\leq T} \triangle M_u \triangle D_u. 
\end{eqnarray}
From Equation \reff{M}, we have
\begin{eqnarray*}
\triangle M_s\neq 0 \,\,\,\mbox{iff} \,\,\, s=T_i=\inf\{u\geq t\,\, \mbox{s.t.} N_u=i\},\,\, i\in \N
\end{eqnarray*}
where $T_i$ is the time of the claim number $i$ and $N$ is the Poisson process representing the number of claims.
If the stopping time $T_i$ is accessible, then from Theorem 15.1 in Rogers and Williams \cite{roge}, we have 
$Z_{T_i^-}=Z_{T_i}$ which is false and so $(T_i)_{i\in \N}$ is a sequence of totally inaccessible stopping times. 
On the other hand, $D$ is predictable process and so from Lemma 27.3 in Rogers and Williams \cite{roge} we have 
$\triangle D_\tau =0$ for every totally inaccessible stopping time $\tau$. 
Since $\triangle M_s=0$ if $s\neq T_i$ and $\triangle D_s=0$ if $s$ is a totally inaccessible stopping time, 
we have 
\begin{eqnarray}\label{var}
\triangle M_s\triangle D_s=0\,\,ds\otimes dP \mbox{ a.s. }
\end{eqnarray}
Since $D$ is nonnegative and nonincreasing, and $Z$, $X$ are
nonnegative, this shows that the process~:
$Y_.X_.- \int_t^. Y_u(\alpha -\beta +(\beta -\int_C\rho_u(z)\,z\, \pi(dz))_{+})du$
is a $P$-local supermartingale, bounded
from below by an $L^1(P)$ random variable, and hence a
$P$-supermartingale. 
\ep
\begin{Lemma} \label{lemelk}
For all $H\in L_+^0(\Fc_T)$, we have~:
\begin{eqnarray} \label{vhegal}
J(H) =
\sup_{Y\in\Yc^0(t)} E\left[Y_TH-
\int_t^T Y_u (\alpha -\beta +(\beta -\int_C\rho_u(z)\,z\, \pi(dz))_{+})du\right]\leq x
\end{eqnarray}
\end{Lemma}
{\bf Proof.}
Fix some $H$ $\in$ $L^0_+(\Fc_T)$. Given  an arbitrary $\tau$ $\in$ $\Tc_t$, 
we define a sequence $(D^n)_{n}$ of elements in $\Dc_t$ by~:
\begin{eqnarray*}
D^n_s\ &=& \exp\left(-\int^{s}_{t}n1_{\tau \leq u}du\right),\;\;
t\leq s\leq T, \; n \in \N.
\end{eqnarray*} 
We then have for all $Z$ $\in$ $\Pc^0(t,x)$~: 
\begin{eqnarray*}
& &E \left[Z_T D^n_T   H  - \int_t^T 
Z_uD^n_u (\alpha -\beta +(\beta -\int_C\rho_u(z)\,z\, \pi(dz))_{+})du \right]\\
& &\leq \sup_{Z\in \Pc^0(t,x),D\in\Dc_t} 
E \left[Z_T D_T   H  - \int_0^T 
Z_uD_u (\alpha -\beta +(\beta -\int_C\rho_u(z)\,z\, \pi(dz))_{+})du \right].
\end{eqnarray*}
Since $D^n_u$ $\rightarrow$ $1_{u\leq \tau}$ a.s., for all $t\leq u\leq T$, 
we have by Fatou's lemma~:
\begin{eqnarray*} 
& & E\left[Z_TH  1_{\tau = T} - \int_0^\tau Z_u (\alpha -\beta +(\beta -\int_C\rho_u(z)\,z\, \pi(dz))_{+})du
\right] \\
&\leq & \sup_{Z\in\Pc^0(t,x),D\in\Dc_t} 
E \left[Z_T D_T   H  - \int_t^T 
Z_uD_u (\alpha -\beta +(\beta -\int_C\rho_u(z)\,z\, \pi(dz))_{+})du \right].
\end{eqnarray*}
Identifying a probability measure $Q$ $\in$ $\Pc^0(t,x)$ with its density
process $Z$, we then obtain from Bayes formula~: 
\begin{eqnarray*} 
J(H) &= & \sup_{Q\in\Pc^0(t,x),\tau\in\Tc_t} 
E^Q\left[H  1_{\tau = T} - \int_0^\tau (\alpha -\beta +(\beta -\int_C\rho_u(z)\,z\, \pi(dz))_{+})du
\right] \\
&\leq & \sup_{Z\in\Pc^0(t,x),D\in\Dc_t} 
E \left[Z_T D_T   H  - \int_0^T 
Z_uD_u (\alpha -\beta +(\beta -\int_C\rho_u(z)\,z\, \pi(dz))_{+})du \right]. 
\end{eqnarray*}
Conversely, by the supermartingale property of 
$Y_.X_.- \int_0^. Y_u(\alpha -\beta +(\beta -\int_C\rho_u(z)\,z\, \pi(dz))_{+})du$
for any $Y$ $\in$ 
$\Yc^0(t)$, see Lemma \ref{lemdc}, we have~:
\begin{eqnarray} \label{interH}
\sup_{Y\in\Yc^0(t)} 
E \left[Z_T D_T   H  - \int_t^T 
Z_uD_u (\alpha -\beta +(\beta -\int_C\rho_u(z)\,z\, \pi(dz))_{+})du \right] \leq x, \; \forall H \in \Cc_+(x).
\end{eqnarray}
Now, by characterization \reff{carac} , any $H$ $\in$ $L_+^0(\Fc_T)$ lies
in $\Cc_+(J(H))$. We then deduce from  inequality \reff{interH} that
\begin{eqnarray*}
\sup_{Y\in\Yc^0(t)} 
E \left[Y_T   H  - \int_0^T 
Y_u (\alpha -\beta +(\beta -\int_C\rho_u(z)\,z\, \pi(dz))_{+})du \right] &\leq& J(H),
\end{eqnarray*}
which proves the required equality \reff{vhegal}.
\ep \\
The dual problem of \reff{pvaleur} is written as:
\begin{eqnarray}\label{valeurm}
\tilde v(t,y) :=\inf_{Y\in\Yc^0(t)} E\left[ \tilde U(yY^{\rho,D}_T) + 
\int_t^T  yY^{\rho,D}_u (\alpha -\beta +(\beta -\int_C\rho_u(z)\,z\,
\pi(dz))_{+})du \right], 
\end{eqnarray}
We shall adopt a dynamic programming principle approach to study the dual value function \reff{valeurm}. 
We recall the dynamic programming principle for our stochastic control problem: for any
stopping time $0\leq \tau\leq T$, $0\leq t\leq T$ and $0\leq h\leq T-t$,
\begin{eqnarray}\label{pppd}
\tilde v(t,y)&=&\Inf_{Y^{\rho,D}\in \Yc^0(t)}
E\left[\tilde v\left((t+h)\wedge \tau,Y^{\rho,D}_{(t+h)\wedge \tau}\right)\right.\\
&+&\left.\int_{t}^{(t+h)\wedge
    \tau}Y^{\rho,D}_u\left(\alpha-\beta+\left(\beta-\int_C\rho_u(z)\,z\,
      \pi(dz)\right)_+\right)du\right],\nonumber 
\end{eqnarray}
where $a\wedge b=\min(a,b)$ ( see e.g. Fleming and Soner \cite{fleson}).\\
\begin{Lemma} Let $t\in [0,T]$ and $Y^{\rho,D}\in \Yc^0(t)$. Then
the process $Y^{\rho,D}$  
evolves according to the following stochastic differential equation
\begin{eqnarray}\label{etat}
dY^{\rho,D}_s&=&Y_{s^-}^{\rho,D}\Big(-dL_s+\int_C(\rho_s(z)-1)\tilde \mu(ds,dz)\Big),
\end{eqnarray}
with 
\begin{eqnarray}\label{controll}
dL_s=-\frac{dD_s}{D_s}1_{\{D_s>0\}},\,\, t\leq s\leq T,\,\,L_{t^-}=0.
\end{eqnarray}
\end{Lemma}
{\bf Proof.}
By It\^o's product rule, we have
\begin{eqnarray*}
dY^{\rho,D}_s
=Z_{s^-}
dD_s
+D_{s^-}dZ_s
+\triangle Z_s\triangle D_s.
\end{eqnarray*}
From Equation \reff{exporetiel}, we have
\begin{eqnarray*}
\triangle Z_s\neq 0\,\,\,\mbox{iff} \,\,\, s=T_i=\inf\{u\geq t\,\, \mbox{s.t.} N_u=i\},\,\, i\in \N.
\end{eqnarray*}
Repeating the same argument as in equation \reff{M},
we have $\triangle Z_s\triangle D_s=0$ $ds\otimes dP \mbox{ a.s. }$ and so equation \reff{etat} is proved.
\ep \\ \\
We denote by $\Lc_t$ the set of adapted processes $(L_s)_{t\leq s\leq T}$ 
with possible jump at time $s=t$ and satisfying equation \reff{controll}.\\
The  Hamilton Jacobi Bellman Variational Inequality arising from the dynamic programming principle 
\reff{pppd} is written as
\begin{eqnarray}\label{HJB2}
&&\min \left \{ 
\frac{\partial \tilde v}{\partial t}(t,y) +
H\left(t,y,\tilde v, \Dy {\tilde v}\right),-\Dy {\tilde v}(t,y)
\right\}=0,\,(t,y)\in [0,T)\times (0,\infty),
\end{eqnarray}
with terminal condition
\begin{eqnarray}\label{HJB2t}
\tilde v(T,y)=\tilde U(y)\,,y\in (0,\infty),
\end{eqnarray} 
where
$$H\left(t,y,\tilde v, \Dy {\tilde v}\right):=\Inf_{\rho \in \Sigma}\left\{A^\rho\left(t,y,\tilde v, \Dy {\tilde v}\right)+
y\left(\alpha-\beta+(\beta-\int_C\rho(z)\,z\, \pi(dz))_+\right)\right\},$$
\begin{eqnarray*}
A^\rho\left(t,y,\tilde v, \Dy {\tilde v}\right) :=
\int_C\left( \tilde v(t,\rho(z) y)- \tilde v(t,y)-(\rho(z) -1)y \Dy {\tilde v}(t,y)\right)\pi(dz), 
\end{eqnarray*} 
 and $\Sigma:=\left\{\rho 
\mbox{ positive Borel function defined on }C 
\mbox{ s.t.} 
 \int_C\Big(|\log\rho(z)| + \rho(z) \Big)\pi(dz) <\infty
\right\}$.
This divides the time-space solvency region $[0,T)\times (0,\infty)$ into a no-jump region
\begin{eqnarray*}
R_1=\left \{
  (t,y)\in [0,T]\times (0,\infty),\mbox{ s.t. }\displaystyle{\frac{\partial \tilde v}{\partial t}(t,y) +
H\left(t,y,\tilde v, \Dy {\tilde v}\right)=0}
\right\}
\end{eqnarray*} 
and a jump region
\begin{eqnarray*}
R_2=\left \{
  (t,y)\in [0,T]\times (0,\infty),\mbox{ s.t. }\displaystyle{\frac{\partial \tilde
      v}{\partial y}(t,y)}=0\right\}.
\end{eqnarray*}

\section{Viscosity solution}

\setcounter{equation}{0}
\setcounter{Assumption}{0}
\setcounter{Example}{0}
\setcounter{Theorem}{0}
\setcounter{Proposition}{0}
\setcounter{Corollary}{0}
\setcounter{Lemma}{0}
\setcounter{Definition}{0}
\setcounter{Remark}{0}
In this section, we provide a rigorous characterization of the dual value 
function $\tilde v$ as a viscosity solution to the HJBVI
(\ref{HJB2}). The function $\tilde v$ is not known to be continuous and 
so we shall work with the notion of discontinuous viscosity solutions. 
We prove that the dual value function lies in the set of functions  
$D_{\gamma}([0,T]\times (0,\infty))$ defined as follows:
\begin{eqnarray*}
D_{\gamma}([0,T]\times (0,\infty)):=\Big\{\displaystyle{f:[0,T]\times (0,\infty)\rightarrow\R}\mbox{ such that },\\
\,\Sup_{y>0} \frac{| f(t,y)|}{y+y^{-\gamma}}<\infty \mbox{ and }\,\Sup_{x>0,y>0} \frac{|f(t,x)-f(t,y)|}{|x-y|(1+x^{-(\gamma+1)}+y^{-(\gamma+1)})}<\infty\Big\}.
\end{eqnarray*}
The following lemma gives some properties of the dual value function $\tilde v$.
\begin{Lemma}\label{classe}
We assume that there exists a solution the the dual problem \reff{valeurm}. 
The following properties hold:\\
1) The dual value function $\tilde v$ is convex in $y$,\\
2) The dual value function $\tilde v$ satisfies the following growth condition
\begin{eqnarray}\label{growthc}
\Sup_{y>0} \frac{|\tilde v(t,y)|}{y+y^{-\gamma}}<\infty
\end{eqnarray}
3)The dual value function $\tilde v$ satisfies 
\begin{eqnarray}\label{acc}
\Sup_{y_1>0,y_2>0} \frac{|\tilde v(t,y_1)-\tilde
  v(t,y_2)|}{|y_1-y_2|(1+y_1^{-(\gamma+1)}+y_2^{-(\gamma+1)})}<\infty.
 \end{eqnarray}
\end{Lemma}
{\bf Proof.}
See Appendix.
\ep \\
\\
\begin{Remark}
When we assume that there exists a solution the the dual problem \reff{valeurm}, 
one can use the conjugate duality relation proved in theorem 5.1 of Mnif and Pham \cite{mnipham01}.
\end{Remark}
Since the dual value 
function $\tilde v$ is locally bounded, the upper and the lower
semi-continuous envelope of the function $\tilde v$ are well-defined. The
definition of  the upper and the lower semi-continuous envelope of a function $\phi$ are given as follows.
\begin{Definition}
(i) The upper semi-continuous envelope of a function $\phi$ is 
\begin{eqnarray}
\phi^*(t,y)= \Limsup_{\stackrel{(t^\prime,y^\prime)\rightarrow (t,y)}{t^\prime\in [0,T),\,y^\prime>0}}
\phi(t^\prime,y^\prime),\,\, \mbox{ for all }(t,y)\in [0,T)\times (0,\infty).
\end{eqnarray} 
(ii) The lower semi-continuous envelope of a function $\phi$  is 
\begin{eqnarray}
\phi_*(t,y)= \Liminf_{\stackrel{(t^\prime,y^\prime)\rightarrow (t,y)}{t^\prime\in [0,T),\,y^\prime>0}}
\phi(t^\prime,y^\prime),\,\, \mbox{ for all }(t,y)\in [0,T)\times (0,\infty).
\end{eqnarray} 
\end{Definition}
Since the continuity of the Hamiltonian $H$ in his arguments is not obvious,
we define lower semi-continuous envelope
of $H$ by \\ $H_*(t,y,\psi,\Dy \psi):=
\Liminf_{\stackrel{(t^\prime,y^\prime)\rightarrow (t,y)}{t^\prime\in [0,T),\,y^\prime>0}}H(t^\prime,y^\prime,\psi,\Dy
\psi)$ for all $(t,y)\in [0,T)\times (0,\infty)$. \\ \\
Adapting the notion of viscosity solutions introduced by Crandall and Lions \cite{cralio86} and then by Soner \cite{son86b}for first integrodifferential operators, we define the viscosity solution as follows:
\begin{Definition}
(i) A function $\phi$ is a viscosity supersolution of (\ref{HJB2}) in  $[0,T)\times (0,\infty)$ if 
\begin{eqnarray}\label{superdef}
&&\min \left\{
\frac{\partial \psi}{\partial t}(\bar t,\bar y) +
H_*\left(\bar t,\bar y,\psi, \Dy \psi \right), -\Dy \psi(\bar t,\bar y) \right\} \leq 0, 
\end{eqnarray} 
whenever $\psi\in C^1([0,T]\times (0,\infty))$ and $\phi_*- \psi$ has a  strict global minimum at 
$(\bar t,\bar y)\in [0,T)\times 
(0,\infty)$.
\\(ii) A function $\phi$ is a viscosity subsolution of (\ref{HJB2}) in  $[0,T)\times (0,\infty)$ if 
\begin{eqnarray}\label{subdef}
&&\min \left\{
\frac{\partial \psi}{\partial t}(\bar t,\bar y) +
H\left(\bar t,\bar y,\psi, \Dy \psi \right), -\Dy \psi(\bar t,\bar y) \right\} \geq 0, 
\end{eqnarray} 
whenever $\psi\in C^1([0,T]\times (0,\infty))$ and $\phi^*-
\psi$ has a  strict global maximum at $(\bar t,\bar y)\in [0,T)\times (0,\infty)$.\\
(iii) A function $\phi$ is a viscosity solution of (\ref{HJB2}) in  $[0,T)\times
(0,\infty)$ if it is both super-solution and sub-solution in  $[0,T)\times (0,\infty)$.
\end{Definition}
\begin{Remark}
As it is seen in the definition of viscosity solutions, we use the lower semi-continuous envelope
of the Hamiltonian H. In fact to prove that the dual value function $\tilde v$ is viscosity solution of our HJBVI, 
we need only the regularity of $H_*$ 
to derive inequality \reff{superbargamma}.
\end{Remark}
The following theorem relates the dual value function $\tilde v$  to the HJBVI \reff{HJB2}. 
\begin{Theorem}\label{viscosite}
The dual value function $\tilde v$ is a viscosity solution of (\ref{HJB2}) in  $[0,T)\times (0,\infty)$.
\end{Theorem}
{\bf Proof.}
See Appendix.
\ep \\
\\
The HJBVI \reff{HJB2} associated to our problem does not provide a complete characterization of the dual value function $\tilde v$. We need to specify the terminal condition. From the definition of $\tilde v$, it's obvious that $\tilde v(T,y)=\tilde U(y)$. Since we use the notion of discontinuous viscosity solutions, we need to characterize  $ \tilde v^*(T,y)$ and $ \tilde v_*(T,y)$ for all $y\in (0,\infty)$ which is the object of the following theorem.
\begin{Theorem}
The terminal conditions of the upper and lower semi-continuous envelope of $\tilde v$ satisfy the following inequalities
\begin{eqnarray}\label{upperterminal}
\tilde v^*(T,y):=\Limsup_{\stackrel{(t^\prime,y^\prime)\rightarrow (T,y)}{t^\prime\in [0,T),\,y^\prime>0}}
\tilde v(t^\prime,y^\prime)\leq \tilde U(y) ,\,\,\mbox{ for all } y\in (0,\infty),
\end{eqnarray}
\begin{eqnarray}\label{lowerterminal}
\tilde v_*(T,y):=\Limsup_{\stackrel{(t^\prime,y^\prime)\rightarrow (T,y)}{t^\prime\in [0,T),\,y^\prime>0}}
\tilde v(t^\prime,y^\prime)\geq \tilde U(y) ,\,\,\mbox{ for all } y\in (0,\infty).
\end{eqnarray}
\end{Theorem}
{\bf Proof.}
We first prove inequality \reff{upperterminal}. Suppose on the contrary that there exists a constant $\eta >0$ such that $\tilde v^*(T,y)\geq \tilde U(y) +2 \eta$. From the definition of $\tilde v^*$, there exists a sequence $((t_n,y_n))_{n\in \N}$ such that $(t_n,y_n)\longrightarrow (T,y)$ and $\tilde v(t_n,y_n)\longrightarrow\tilde v^*(T,y)$ when $n$ tends to infinity, which implies
\begin{eqnarray*}
E[\tilde
U(y_nY_{T}^{\rho_n,D_n})+\int_{t_n}^{T}y_nY_s^{\rho_n,D_n}(\alpha-\beta+(\beta-\int_C\rho_{n_{s}}(z)\,z\,\pi(dz))_+)ds]\geq
\tilde U(y) +\eta,
\end{eqnarray*}
for all $Y^{\rho_n,D_n}\in \Yc^0(t_n)$.Choosing $\rho_{n_{s}}=1$ and $D_{n_{s}}=1$ for all $s\in[t_n,T]$, we obtain
\begin{eqnarray*}
\tilde U(y_n)+\int_{t_n}^{T}y_n(\alpha-\beta+(\beta-\int_C\,z\,\pi(dz))_+)ds\geq \tilde U(y) +\eta. 
\end{eqnarray*}
Sending $n$ to infinity, we have $\tilde U(y) \geq \tilde U(y) +\eta$ which is wrong and so inequality \reff{upperterminal} is proved.\\
We prove now inequality \reff{lowerterminal}. From the definition of $\tilde v_*$, there exists a sequence $((t_n,y_n))_{n\in \N}$ such that $(t_n,y_n)\longrightarrow (T,y)$ and $\tilde v(t_n,y_n)\longrightarrow\tilde v_*(T,y)$ when $n$ tends to infinity. \\
From the definition of the dynamic version of the value function, we have for all $x\geq (\beta-\alpha)(T-t_n)$
\begin{eqnarray}\label{dyntn}
v(t_n,x)\geq U(x+(\alpha-\beta)(T-t_n)).
\end{eqnarray} 
Let $\epsilon>0$, there exists $n_0\in \N$ such that for all $n\geq n_0$, we have 
\begin{eqnarray}\label{dyntnn}
 U(x+(\alpha-\beta)(T-t_n))\geq U(x)-\epsilon.
\end{eqnarray} 
Using the conjugate duality relation of Theorem 5.1 of Mnif and Pham \cite{mnipham01} and relations \reff{dyntn} and \reff{dyntnn}, we obtain
\begin{eqnarray}\label{dyntnnn} 
\tilde v(t_n,y_n)&=&\Max_{x\geq (\beta-\alpha)(T-t_n) }[v(t_n,x)-xy_n]\nonumber\\ &\geq& \Max_{x\geq (\beta-\alpha)(T-t_n) }[U(x)-xy_n]-\epsilon.
\end{eqnarray}
Since $(t_n,y_n)\longrightarrow (T,y)$ when $n$ tends to infinity, there exists $n_1\in \N$ such that for all $n\geq n_1$, $I(y_n)\geq (\beta-\alpha)(T-t_n)$ and so $\Max_{x\geq (\beta-\alpha)(T-t_n) }[U(x)-xy_n]=\tilde U(y_n)$.
For $n\geq n_0 \vee n_1$, inequality \reff{dyntnnn} implies 
\begin{eqnarray*}
\tilde v(t_n,y_n)\geq \tilde U(y_n)-\epsilon.
\end{eqnarray*}
Sending $n$ to infinity and $\epsilon$ to $0$, we prove inequality\reff{lowerterminal}.
\ep 
\section{Uniqueness}

\setcounter{equation}{0}
\setcounter{Assumption}{0}
\setcounter{Example}{0}
\setcounter{Theorem}{0}
\setcounter{Proposition}{0}
\setcounter{Corollary}{0}
\setcounter{Lemma}{0}
\setcounter{Definition}{0}
\setcounter{Remark}{0}
We turn now to uniqueness questions. Our next main result is a comparison
principle for discontinuous viscosity solutions to the HJBVI \reff{HJB2}. 
It states that we can compare a viscosity sub-solution and a viscosity 
super-solution to the HJBVI \reff{HJB2} on $[0,T)\times (0,\infty)$, provided that we can compare them
at terminal date as usual in parabolic problems. \\
The main difficulty in the comparison theorem comes from the discontinuity
of the  Hamiltonian. 
Here, we assume that the claims take values in the set 
$C=\{\delta_1,\delta_2,...,\delta_d\}$, $\delta_i>0$, $1\leq i\leq d$. In this case the Hamiltonian contains an inf on a bounded set which makes it continuous.
We denote by $\pi_i$, $1\leq i\leq d$ the intensity of the Poisson process associated 
to the claim having the size $\delta_i$. The set $\Sigma$ is defined as follows~:
\begin{eqnarray}\label{sigmaprime}
\Sigma:=\left\{ \rho =(\rho_i)_{1\leq i\leq d},\,\, \rho_i>0,\,1\leq i\leq d\right\}.
\end{eqnarray}
The Hamiltonian $H$ is given by
\begin{eqnarray}\label{tildeH}
H\left(t,y,\tilde v, \Dy {\tilde v}\right)=\Inf_{\rho \in \Sigma}\left\{A^\rho\left(t,y,\tilde v, \Dy {\tilde v}\right)+
y\left(\alpha-\beta+(\beta-\Sum_{i=1}^d\rho_i \delta_i \pi_i)_+\right)\right\},
\end{eqnarray} 
where 
\begin{eqnarray*}
A^\rho\left(t,y,\tilde v, \Dy {\tilde v}\right) :=\Sum_{i=1}^d \pi_i\left( \tilde v(t,\rho_i y)- 
\tilde v(t,y)-(\rho_i -1)y \Dy {\tilde v}(t,y)\right).
\end{eqnarray*} 
The following lemma sates some properties of the functions  $\tilde v_*$ and  $\tilde v^*$.
\begin{Lemma}\label{conv}
We assume that there exists a solution the the dual problem \reff{valeurm}. 
The following properties hold~:\\
1) The functions $\tilde v^*$ and $\tilde v_*$ are convex in $y$.\\
2) The functions $\tilde v^*$ and $\tilde v_*$ are nonincreasing on $(0,\infty)$.\\
3) The functions $\tilde v_*$ and  $\tilde v^*$ satisfy the following growth condition
\begin{eqnarray}\label{growthc}
\Sup_{y>0} \frac{|\tilde v_*(t,y)|}{y+y^{-\gamma}},\,\,\Sup_{y>0} \frac{|\tilde v^*(t,y)|}{y+y^{-\gamma}}<\infty
\end{eqnarray}
4) The functions $\tilde v_*$ and  $\tilde v^*$ satisfy
\begin{eqnarray*}
\Sup_{y_1>0,y_2>0} 
\frac{|\tilde v_*(t,y_1)-\tilde v_*(t,y_2)|}
{|y_1-y_2|(1+y_{1}^{-(\gamma+1)}+y_{2}^{-(\gamma+1)})}<\infty.
\end{eqnarray*}
\begin{eqnarray*}
\Sup_{y_1>0,y_2>0} 
\frac{|\tilde v^*(t,y_1)-\tilde v^*(t,y_2)|}
{|y_1-y_2|(1+y_{1}^{-(\gamma+1)}+y_{2}^{-(\gamma+1)})}<\infty.
 \end{eqnarray*}
\end{Lemma}
{\bf Proof.}
1) Fix $\lambda \in (0,1)$ and $(y,y^{'},y^{''})\in (0,\infty)^3$ such that $y=\lambda y^{'}+(1-\lambda)y^{''}$. 
From the definition of $\tilde v^*(t,y)$, there exists a sequence $\big((t_n,y_n)\big)_n$ such that 
$v(t_n,y_n)\longrightarrow \tilde v^*(t,y)$ when $n$ goes to infinity. 
We set $y_n=\lambda y^{'}+(1-\lambda)y^{''}_n$, then we have 
$y^{''}_n\longrightarrow y^{''}$ when $n$ goes to infinity.
Since $\tilde v$ is convex in $y$, then we have
\begin{eqnarray}\label{v*conv}
\tilde v(t_n,y_n)&\leq& \lambda\tilde v(t_n,y^{'})+(1-\lambda)\tilde v(t_n,y^{''}_n)
\end{eqnarray}
Sending $n$ to infinity, inequality \reff{v*conv} implies
\begin{eqnarray}
\tilde v^*(t,y)\leq  \lambda\tilde v^*(t,y^{'})+(1-\lambda)\tilde v^*(t,y^{''}),
\end{eqnarray} 
which is the desired result.\\
We turn to the convexity of $\tilde v_*$.
From Theorem 5.1 of Mnif and Pham \cite{mnipham01}, we have 
\begin{eqnarray}\label{conv1}
\tilde v(t,y)=\Max_{x\geq (\beta-\alpha)(T-t)}[v(t,x)-xy].
\end{eqnarray} 
We fix $(t,y)\in [0,T]\times (0,\infty)$. From the definition of $\tilde v_*$, 
there exists a sequence$\big((t_n,y_n)\big)_n$ such that 
$v(t_n,y_n)\longrightarrow \tilde v_*(t,y)$ and $(t_n,y_n)\longrightarrow (t,y)$ when $n$ goes to infinity. 
From equation \reff{conv1}, we have $\tilde v(t_n,y_n)\geq v(t_n,x)-xy_n \geq v_*(t,x)-xy_n $ 
for all $x\geq (\beta-\alpha)(T-t_n)$. Sending $n$ to $\infty$, 
we have $\tilde v_*(t,y)\geq v_*(t,x)-xy$ for all $x\geq (\beta-\alpha)(T-t)$ and so
\begin{eqnarray}\label{conv2}
\tilde v_*(t,y)\geq \Max_{x\geq (\beta-\alpha)(T-t)}[v_*(t,x)-xy].
\end{eqnarray} 
From Theorem 5.1 of Mnif and Pham \cite{mnipham01}, there exists $\hat x$ such that 
\begin{eqnarray}\label{conv3}
\tilde v(t,y)= v(t,\hat x)-\hat xy
\end{eqnarray} 
 and $\Dy{\tilde v}(t,y)=-\hat x$.
We consider a sequence $\big((t_n,y_n)\big)_n$ such that 
$v(t_n,y_n)\longrightarrow \tilde v_*(t,y)$ and $(t_n,y_n)\longrightarrow (t,y)$ when $n$ goes to infinity.
From equation \reff{conv3}, there exists $\hat x_n$ such that $\tilde v(t_n,y_n)= v(t_n,\hat x_n)-\hat x_ny_n$. Since 
$\tilde v \in D_{\gamma}([0,T]\times (0,\infty))$, we have 
$\hat x_n:=|\Dy{\tilde v}(t_n,y_n)|\leq K(1+y_n^{-(\gamma+1))})$ where $K$ is a 
positive constant independent of $n$. Since $y_n\longrightarrow y$ when $n$ goes to infinity, 
the sequence $(\hat x_n)_n$ is bounded, and so long a subsequence 
$\hat x_n\longrightarrow \hat x $ when $n$ goes to infinity and so we have
\begin{eqnarray}\label{conv4}
\Lim_{n\longrightarrow \infty}v(t_n,\hat x_n)= \tilde v_*(t,y)+\hat xy= v_*(t,\hat x)
\end{eqnarray} 
From inequality \reff{conv2} and equation \reff{conv4}, we deduce that
\begin{eqnarray*}
\tilde v_*(t,y)= \Max_{x\geq (\beta-\alpha)(T-t)}[v_*(t,x)-xy].
\end{eqnarray*}
and so $\tilde v_*$ is convex in $y$.\\
2) Fix $(t,y)\in [0,T)\times (0,\infty)$. Since at time $t$, only $D\in \Dc_t$ could make a jump, 
we have $y=Y_t\geq Y_{t^+}$. From the dynamic programming principle
\reff{pppd} we have
\begin{eqnarray*}
\tilde v(t,y)\leq \tilde v(t,y(1-\delta))\,\mbox{ for all }0<\delta<1,
\end{eqnarray*}
and so the dual value
function $\tilde v$ in non-increasing with respect to $y$. 
This yields that $\tilde v_*$ and $\tilde v^*$ are nonincreasing.
\\
3) Fix $(t,y)\in [0,T)\times (0,\infty)$. From the definition of $\tilde v^*(t,y)$, there exists 
a sequence $\big((t_n,y_n)\big)_n$ such that 
$v(t_n,y_n)\longrightarrow \tilde v_*(t,y)$ when $n$ goes to infinity. From the growth condition of $\tilde v$, we have
\begin{eqnarray*}
|\tilde v(t_n,y_n)| \leq K(y_n+y_n^{-\gamma}),
\end{eqnarray*} 
where $K$ is a positive constant. Sending $n$ to infinity, we obtain the desired result. We use similar arguments to prove that $\Sup_{y>0} \frac{|\tilde v^*(t,y)|}{y+y^{-\gamma}}<\infty$.\\
4) We prove only the first inequality. The second one is obtained by using similar arguments.
Fix $(t,y_1,y_2)\in [0,T)\times (0,\infty)\times (0,\infty)$. From the definition of $\tilde v^*(t,y_2)$, there exists 
a sequence $\big((t_n,y_{2,n})\big)_n$ such that 
$v(t_n,y_{2,n})\longrightarrow \tilde v_*(t,y_2)$ when $n$ goes to infinity. From inequality \reff{acc}, we have
\begin{eqnarray*}
\tilde v(t_n,y_{1})-\tilde
  v(t_n,y_{2,n})\leq K |y_{1}-y_{2,n}|(1+y_{1}^{-(\gamma+1)}+y_{2,n}^{-(\gamma+1)}), 
\end{eqnarray*} 
where $K$ is a positive constant.
Since $ \Liminf_{n\longrightarrow \infty}\tilde v(t_n,y_{1})\geq\tilde v_*(t,y_{1}) $, 
we obtain after sending $n$ to infinity
\begin{eqnarray*}
-\tilde v_*(t,y_{2}) \leq K|y_{1}-y_{2}|(1+y_{1}^{-(\gamma+1)}+y_{2}^{-(\gamma+1)}) - \tilde v_*(t,y_{1}). 
\end{eqnarray*} 
Using similar arguments, we deduce the inverse inequality and so
\begin{eqnarray*}
\Sup_{y_1>0,y_2>0} 
\frac{|\tilde v_*(t,y_1)-\tilde v_*(t,y_2)|}
{|y_1-y_2|(1+y_{1}^{-(\gamma+1)}+y_{2}^{-(\gamma+1)})} <\infty.
 \end{eqnarray*}
\ep \\
\begin{Remark}
{\rm One could prove the monotinicity of $\tilde v^* $ by using viscosity solutions arguments. In fact
for each $\epsilon >0$, we define $W(t,y)=\tilde v(t,y)-\epsilon y$, 
$(t,y)\in [0,T)\times (0,\infty)$. The function $W$ satisfies in the viscosity sense 
$\Dy {W} \leq -\epsilon$, i.e. for all $(y_0,\psi)\in ((0,\infty),C^1([0,T]\times (0,\infty))) $ such that
\begin{eqnarray}\label{devr}
(W^*-\psi)(t,y_0)=\Max_{y\in (0,\infty)}(W^*-\psi)(t,y)
\end{eqnarray} 
and so we have $\Dy {\psi}(t,y_0)\leq -\epsilon$. This proves that $\psi(t,.)$ is nonincreasing on a neighborhood $V(y_0)$.
Let $(y_1,y_2)\in V(y_0)$, we want to prove that $W^*(t,y_1)\geq W^*(t,y_2)$.\\
Suppose that $W^*(t,y_1)<W^*(t,y_2)$. We consider the function $V(t,y)=W^*(t,y_1)$ which solves 
\begin{eqnarray}\label{equaaux}
\Dy {V}=0 \mbox{ on }(y_1,y_2),
\end{eqnarray} 
together with the boundary conditions $V(t,y_1)=V(t,y_2)=W^*(t,y_1)$. Since $W$ is a viscosity subsolution of 
Equation \reff{equaaux}, From the comparison theorem ( see Barles \cite{B}, Theorem 2.7 in the case of continuous viscosity solutions), we have
\begin{eqnarray*}
\Inf_{[y_1,y_2]}(V-W^*)(t,y)=\Inf((V-W^*)(t,y_1),(V-W^*)(t,y_2))=0,
\end{eqnarray*} 
which implies $V(t,y)\geq W^*(t,y)$ for all $y\in [y_1,y_2]$. 
Since $W^*(t,y_0)\leq V(t,y_0)=W^*(t,y_1)$ and  $\psi(t,y_0)>\psi(t,y_1)$, then 
\begin{eqnarray*}
(W^*-\psi)(t,y_0)<(W^*-\psi)(t,y_1),
\end{eqnarray*} 
which contradicts \reff{devr}. Sending $\epsilon \longrightarrow 0^+$, we obtain the desired result.\\}
\end{Remark}
Now, we are able to prove  the following comparison principle~: 
\begin{Theorem}\label{comparaison}
Let $\tilde v_1$ (resp $\tilde v_2$) $\in D_\gamma ([0,T]\times (0,\infty))$ be a viscosity subsolution (resp supersolution) of \reff{HJB2} in $[0,T]\times (0,\infty)$  such that $\tilde v_1^*(T,y)\leq \tilde v_{2*}(T,y)$. We assume that 
$\tilde v_{2*}$ is convex and nonincreasing in $y$, then
\begin{eqnarray}\label{comp}
\tilde v_1^*(t,y)\leq \tilde v_{2*}(t,y),\,\mbox{ for all }(t,y)\in [0,T]\times (0,\infty)
\end{eqnarray}
\end{Theorem}
{\bf Proof.}
See Appendix.
\ep \\ \\
By combining the previous results, we finally obtain the following PDE  characterization of the dual value function.
\begin{Corollary}\label{coro}We assume that there exists a solution the the dual problem \reff{valeurm}. 
The dual value function $\tilde v$ is the unique viscosity solution  of \reff{HJB2} with terminal condition $\tilde v(T,y)=\tilde U(y)$ in the class of functions $D_\gamma ([0,T]\times (0,\infty))$.
\end{Corollary}
\begin{Remark}
In this remark, we formally discuss the numerical implications
of Corollary \ref{coro}. We can solve numerically the associated HJBVI by using an algorithm based on policy iterations. 
Then thanks to a verification theorem, We characterize the optimal insurance strategy by the solution 
of the variational inequality. These results are the object of the paper Mnif\cite{mnif10}.
\end{Remark}
\section{Appendix}

\setcounter{equation}{0}
\setcounter{Assumption}{0}
\setcounter{Example}{0}
\setcounter{Theorem}{0}
\setcounter{Proposition}{0}
\setcounter{Corollary}{0}
\setcounter{Lemma}{0}
\setcounter{Definition}{0}
\setcounter{Remark}{0}
\subsection{Proof of Lemma \ref{classe}}
$1)$ From Theorem 5.1 of Mnif and Pham \cite{mnipham01}, we have 
$\tilde v(t,y)=\Max_{x\geq (\beta-\alpha)(T-t)}[v(t,x)-xy]$. The  convexity property of $\tilde v$ in $y$ holds since 
it is the upper envelope of affine functions.\\
$2)$ Since the controls $\rho_s=1$ and $D_s=1$, $s\in [t,T]$ lie in $\mathcal{U}_t \times \mathcal{D}_{t}$, we have
\begin{eqnarray}\label{minorant4}
\tilde v(t,y)\leq \tilde U(y)+Ky, 
\end{eqnarray}
where $K$ is a constant.\\
Let $(Z^n:=Z^{\rho^n},D^n)$ be a minimizing sequence of $\tilde v(t,y)$.
From the definition of these minimizing sequences,
there exist $\epsilon_n$ and $n_0\in \N$ such that
$\epsilon_n \longrightarrow 0$ when  $n \longrightarrow \infty$ and for all $n\geq n_0$, we have
\begin{eqnarray}\label{minorant}
\tilde v(t,y) &\geq & E\left[ \tilde U(yZ^n_TD^n_T)\right] \nonumber\\
&+& y E\left[\int_t^T  Z^n_uD^n_u (\alpha -\beta +(\beta -\int_C\rho^n_u(z)\,z\, \pi(dz))_{+})du \right]-\epsilon_n.
\end{eqnarray}
Since  $\epsilon_n \longrightarrow 0$ when  $n\longrightarrow \infty$, there exists $n_1\in \N$ such
that for all $n\geq n_1$, we have $\epsilon_n\leq \tilde U(y)+y$. We recall
That $\tilde U(y)\geq U(0^+)\geq 0$ and so $\tilde U(y)+y>0$ since $y>0$.
Using the boundedness of $D^n$, Jensen's inequality and the martingale property of $Z^n$, we have:
\begin{eqnarray}\label{minorant3}
E\left[\tilde U(yZ^n_TD^n_T)\right]&\geq& \tilde U(yE\left[Z^n_T\right])\nonumber\\
&\geq& \tilde U(y).
\end{eqnarray}
For the second term of the r.h.s of
inequality \reff{minorant}, since $D^n_s\leq 1$ for all $s\in [t,T]$, using
Fubini's theorem and the martingale property of $Z^n$, we have
\begin{eqnarray}\label{minorant2}
E\left[\int_t^T y  Z^n_uD^n_u  (\alpha -\beta +(\beta
  -\int_C\rho_u(z)\,z\, \pi(dz))_{+})du \right]&\geq& y(\alpha-\beta)E\left[ \int_t^T   Z^n_uD^n_u du\right]\nonumber\\
&\geq& y(\alpha-\beta) \int_t^T E\left[Z^n_u\right] du\nonumber\\
&\geq& K^\prime y,
\end{eqnarray}
where  $K^\prime$ is a constant independent of $y$. Inequalities \reff{minorant3} and \reff{minorant2} imply that
\begin{eqnarray}\label{minorant5}
\tilde v(t,y)\geq  \tilde U(y)+ K^{\prime}y.
\end{eqnarray}
From inequalities \reff{minorant4} and \reff{minorant5}, we deduce that 
\begin{eqnarray}\label{comportement}
\Sup_{y>0} \frac{|\tilde v(t,y)|}{y+\frac{y^{-\gamma}}{\gamma}}<\infty
\end{eqnarray}
$3)$ Using the convexity in $y$ of the dual value function $\tilde v$, we have
\begin{eqnarray}\label{accroissement}
\frac{|\tilde v(t,y_1)-\tilde v(t,y_2)|}{|y_1-y_2|}\leq |\tilde v^{\prime}_d(t,y_1)|+|\tilde v^{\prime}_d(t,y_2)|,
\end{eqnarray}
where $\tilde v^{\prime}_d$ is the right-hand derivative with respect to
the variable $y$. Let $(Z^n:=Z^{\rho^n},D^n)$ be a minimizing sequence of $\tilde v(t,y_1)$.  Let
$\delta >0$ and $(Z^{\prime n}:=Z^{\prime \rho^n},D^{\prime n})$ be a minimizing sequence of
$\tilde v(t,y_1+\delta)$. From the definition of these minimizing sequences,
there exist $\epsilon_n$, $\epsilon^{\prime}_n$ and $n_0\in \N$ such that
$\epsilon_n \longrightarrow 0$, $\epsilon^{\prime}_n \longrightarrow 0$ when  $n
\longrightarrow \infty$ and for all $n\geq n_0$, we have
\begin{eqnarray}\label{accroissement1}
\tilde v(t,y_1) &\geq&  y_1^{-\gamma}E\left[ \frac{(Z^n_TD^n_T)^{-\gamma}}{\gamma}\right] \nonumber\\
&+& y_1 E\left[\int_t^T  Z^n_uD^n_u 
\Big(\alpha -\beta +(\beta -\int_C\rho^n_u(z)\,z\, \pi(dz))_{+}\Big)du \right]-\epsilon_n,
\end{eqnarray}
and
 \begin{eqnarray}\label{accroissement5}
\tilde v(t,y_1+\delta)
 &\geq&  (y_1+\delta)^{-\gamma}E\left[ \frac{(Z^{\prime n}_TD^{\prime
    n}_T)^{-\gamma}}{\gamma}\right] \nonumber\\
&+& (y_1+\delta) E\left[\int_t^T  Z^{\prime n}_uD^{\prime n}_u
   \Big (\alpha -\beta +(\beta -\int_C\rho^{\prime n}_u(z)\,z\, \pi(dz))_{+}\Big)du
    \right]\nonumber\\
&-&\epsilon^\prime _n.
\end{eqnarray}
Using the definition of $\tilde v(t,y_1)$ and $\tilde v(t,y_1+\delta)$, we have
\begin{eqnarray*}
\tilde v(t,y_1) \leq  y_1^{-\gamma}E\left[\frac{( Z^{\prime n}_TD^{\prime
    n}_T)^{-\gamma}}{\gamma}\right] + y_1 E\left[\int_t^T  Z^{\prime n}_uD^{\prime n}_u
   \Big (\alpha -\beta +(\beta -\int_C\rho^{\prime n}_u(z)\,z\, \pi(dz))_{+}\Big)du \right],
\end{eqnarray*}
and
\begin{eqnarray}\label{accroissement2}
\tilde v(t,y_1+\delta)& \leq & (y_1+\delta)^{-\gamma}E\left[\frac{( Z^{
    n}_TD^{n}_T)^{-\gamma}}{\gamma}\right] \\
&+& (y_1+\delta) E\left[\int_t^T  Z^{n}_uD^{n}_u
    \Big(\alpha -\beta +(\beta -\int_C\rho^{n}_u(z)\,z\, \pi(dz))_{+}\Big)du \right].\nonumber
\end{eqnarray}
Using inequalities \reff{accroissement1} and \reff{accroissement2}, we deduce
that
\begin{eqnarray}\label{derivedte}
& &\frac{\tilde v(t,y_1+\delta)-\tilde v(t,y_1)}{\delta}\nonumber\\
&\leq&\frac{(y_1+\delta)^{-\gamma}-y_1^{-\gamma}}{\delta}E\left[\frac{(
  Z^{n}_TD^{n}_T)^{-\gamma}}{\gamma}\right]\nonumber\\
&+& E\left[\int_t^T  Z^{n}_uD^{n}_u
    \Big(\alpha -\beta +(\beta -\int_C\rho^{n}_u(z)\,z\, \pi(dz))_{+}\Big)du \right]
+\frac{\epsilon_n}{\delta}.
\end{eqnarray}
From inequality \reff{minorant3}, we have
\begin{eqnarray}\label{minorant11}
E\left[(Z^{n}_TD^{n}_T)^{-\gamma}\right]\geq K, 
\end{eqnarray}
where $K$ is a positive
constant independent of $y$ and $\delta$. Since  $\epsilon_n
\longrightarrow 0$ when  $n\longrightarrow \infty$, there exists $n_1$ such
that for all $n\geq n_1$, we have $\frac{\epsilon_n}{\delta}\leq 1$.
For the second term of the r.h.s of
inequality \reff{derivedte}, since $D_s\leq 1$ for all $s\in [t,T]$, using Fubini and martingale property of $Z$, we have
\begin{eqnarray}\label{minorant10}
E\left[\int_t^T   Z^n_uD^n_u  \Big(\alpha -\beta +(\beta
  -\int_C\rho^n_u(z)\,z\, \pi(dz))_{+}\Big)du \right]
&\leq& \alpha E\left[ \int_t^T   Z^n_uD^n_u du\right]\nonumber\\
&\leq& \alpha \int_t^T E\left[Z^n_u\right] du\nonumber\\
&\leq& K .
\end{eqnarray}
Using  inequalities \reff{derivedte}, \reff{minorant11} and \reff{minorant10}
\begin{eqnarray*}
\frac{\tilde v(t,y_1+\delta)-\tilde v(t,y_1)}{\delta} \leq K(\frac{(y_1+\delta)^{-\gamma}-y_1^{-\gamma}}{\delta}+1).
\end{eqnarray*}
Sending $\delta \longrightarrow 0^+$, we have
\begin{eqnarray}\label{derivedte2}
\tilde v^{\prime}_d(t,y_1) \leq K(-y_1^{-(\gamma+1)}+1).
\end{eqnarray}
Similarly, we obtain that
\begin{eqnarray}\label{derivedte3}
& &\frac{\tilde v(t,y_1+\delta)-\tilde v(t,y_1)}{\delta}\nonumber\\
&\geq&\frac{(y_1+\delta)^{-\gamma}-y_1^{-\gamma}}{\delta}E\left[\frac{(
  Z^{\prime n}_TD^{\prime n}_T)^{-\gamma}}{\gamma}\right]\nonumber\\
&+& E\left[\int_t^T  Z^{\prime n}_uD^{\prime n}_u
   \Big (\alpha -\beta +(\beta -\int_C\rho^{\prime n}_u(z)\,z\, \pi(dz))_{+}\Big)du \right]-\frac{\epsilon^\prime_n}{\delta}.
\end{eqnarray}
We define $\tilde \rho_s:=C^{\tilde \rho}=\frac{\beta}{\int_Cz\pi(dz)}$ and $\tilde D_s=1$ for
all $s\in[t,T]$. 
From the definition of $\tilde v$ and using the martingale property of $Z^{\tilde
  \rho}$ , we have
\begin{eqnarray}\label{acc1}
 & &\tilde v(t,y_1+\delta)\nonumber\\ &\leq&  
(y_1+\delta)^{-\gamma} 
E\left[\frac{(
  Z^{\tilde \rho}_T \tilde D_T)^{-\gamma}}{\gamma}\right]+ 
 (y_1+\delta) E\left[\int_t^T   Z^{\tilde \rho}_u \tilde D_u
   \Big (\alpha -\beta +(\beta -\int_C\tilde \rho_u(z)\,z\, \pi(dz))_{+}\Big)du \right]\nonumber\\
&=&(y_1+\delta)^{-\gamma} 
E\left[\frac{(
  Z^{\tilde \rho}_T )^{-\gamma}}{\gamma}\right]+ 
(y_1+\delta) (\alpha -\beta ) (T-t).
\end{eqnarray}
From inequality \reff{accroissement5}, using the martingale of $Z^{\prime n}$ and 
since $0 \leq D^{\prime n}_s\leq 1$ for all $s\in [t,T]$,  we have
\begin{eqnarray}\label{acc2}
\tilde v(t,y_1+\delta) &\geq& (y_1+\delta)^{-\gamma}E\left[(\frac{Z^{\prime n}_TD^{\prime n}_T}{\gamma})^{-\gamma}\right]
+  (y_1+\delta) (\alpha -\beta ) (T-t)-\epsilon^\prime_n
\end{eqnarray}
from \reff{acc1} and \reff{acc2}, we deduce that
\begin{eqnarray*}
E\left[(Z^{\prime n}_TD^{\prime n}_T)^{-\gamma}\right]\leq E\left[(Z^{\tilde
    \rho}_T )^{-\gamma}\right]+\gamma \epsilon^\prime_n (y_1+\delta)^\gamma.
\end{eqnarray*}
We know that $Z^{\tilde \rho}$ is given by the formula \footnote{ The solution of the SDE 
$dZ_t=Z_{t^-}dH_t$ is given by the Dol\'eans Dade exponential formula 
$Z_t=\exp{(H_t-\frac{1}{2}<H^c>_t)}{\displaystyle \prod_{0\leq s\leq t}}(1+\triangle H_s)\exp{(-\triangle H_s)}
$ where $<H^c>$ is the quadratic variation of the continuous part of $H$.} 
\begin{eqnarray*}
Z^{\tilde \rho}_T=\exp(\varrho(t-T)(C^{\tilde \rho}-1))(C^{\tilde \rho})^{N_T-N_t}
\end{eqnarray*}
and so $ E\left[(Z^{\tilde
    \rho}_T )^{-\gamma}\right]<\infty$.
Since  $\epsilon_n^\prime \longrightarrow 0$ when  $n\longrightarrow \infty$,
there exists $n_2$ such that for all $n\geq n_2$, we have $E\left[(Z^{\prime
    n}_TD^{\prime n}_T)^{-\gamma}\right]\leq K^\prime$, where $K^\prime$ is a
positive constant independent of $y_1$ and $\delta$. 
Using the boundedness of $D$, inequality \reff{derivedte3} implies
\begin{eqnarray*}
\frac{\tilde v(t,y_1+\delta)-\tilde v(t,y_1)}{\delta} 
\geq K^\prime\Big(\frac{(y_1+\delta)^{-\gamma}-y_1^{-\gamma}}{\delta}-1\Big).
\end{eqnarray*}
Sending $\delta \longrightarrow 0^+$, we have
\begin{eqnarray}\label{derivedte4}
\tilde v^{\prime}_d(t,y_1) \geq -K^\prime(y_1^{-(\gamma+1)}+1).
\end{eqnarray}
Using \reff{accroissement}, \reff{derivedte2} and \reff{derivedte4}, we obtain
\begin{eqnarray*}
\Sup_{y_1>0,y_2>0} \frac{|\tilde v(t,y_1)-\tilde
  v(t,y_2)|}{|y_1-y_2|\Big(1+y_1^{-(\gamma+1)}+y_2^{-(\gamma+1)}\Big)}<\infty
\end{eqnarray*}
\ep 

\subsection{Proof of Theorem \ref{viscosite}}
We first prove that  $\tilde v$  is a viscosity sub-solution of
(\ref{HJB2}) in $[0,T)\times (0,\infty)$. 
\\Let  $(t,y)\in [0,T)\times (0,\infty)$ and $\psi\in C^1([0,T]\times (0,\infty))$
such that without loss of generality 
\begin{eqnarray*}
0=(\tilde v^*-\psi)(t,y)=\Max_{[0,T)\times (0,\infty)}(\tilde v^*-\psi).
\end{eqnarray*}
From the definition of $\tilde v^*$, there exists a sequence $(t_n,y_n)\in
[0,T)\times (0,\infty) $  such that \\ $(t_n,y_n)\longrightarrow (t,y)$ and $\tilde
v(t_n,y_n)\longrightarrow \tilde 
v^*(t,y)$ when $n\longrightarrow \infty$.\\ 
For  $\eta >0$, ${\rho_n}_s=\tilde \rho $ a positive Borel function and ${D_{n}}_s=1$ for all $s\geq t_n$, we set 
\begin{eqnarray*} 
\theta_n:=\inf\{t\geq t_n\mbox{ such that }(t,y_nY_t^{\rho_n,D_n})\notin B_\eta(t_n,y_n)\}\wedge T,
\end{eqnarray*} 
where $B_{\eta}(t_n, y_n)=\{(t,y)\in [0,T]\times (0,\infty) \mbox{ such that }|t-t_n|+|y-y_n|\leq \eta
\}$. By the right continuity of the paths, we have $\theta_n>t_n$ a.s. 
For all $0<h<T-t_n$, the dynamic programming principle
\begin{eqnarray*}
\tilde v(t_n,y_n)&=&\Inf_{Y^{\rho,D}\in \Yc^0(t_n)}E\left[\tilde
  v\left((t_n+h)\wedge \theta_n, y_n Y^{\rho,D}_{(t_n+h)\wedge
      \theta_n}\right)\right.\\ 
&+&\left.\int_{t_n}^{(t_n+h)\wedge
    \theta_n} y_n Y^{\rho,D}_s\left(\alpha-\beta+\left(\beta-\int_C\rho_s(z)\,z\,
      \pi(dz)\right)_+\right)ds\right] 
\end{eqnarray*} 
implies
\begin{eqnarray}\label{ppdgam} 
\gamma_n +\psi(t_n,y_n)&\leq& E\left[\psi \left((t_n+h)\wedge
    \theta_n, y_n Y^{\rho_n, D_n}_{(t_n+h)\wedge \theta_n}\right)\right.\\ 
&+&\left.\int_{t_n}^{(t_n+h)\wedge \theta_n} y_n Y^{\rho_n,D_n}_s\left(\alpha-\beta
    +\left(\beta-\int_C\rho_n(z)\,z\, \pi(dz)\right)_+\right)ds\right]\nonumber, 
\end{eqnarray}
where the sequence $\gamma_n:= \tilde v(t_n,y_n)-\psi(t_n,y_n)$ is determinist
and converges to zero when $n$ tends to infinity. 
Applying It\^o's formula to $\psi (t_n+h,y_nY^{\rho_n,D_n}_{t+h})$, we get
\begin{eqnarray*}
& &E\left[\frac{1}{h}\int_{t_n}^{(t_n+h)\wedge \theta_n}\frac{\partial
    \psi}{\partial t}\left(s, y_n Y^{\rho_n,D_n}_s\right) +
  A^{\rho_{n}}\left(s, y_n Y^{\rho_n,D_n}_s,\psi, \Dy \psi\right)ds\right]\\ 
&+&E\left[
\frac{1}{h}\int_{t_n}^{(t_n+h)\wedge
    \theta_n} y_nY^{\rho_n,D_n}_s\left(\alpha-\beta+\left(\beta-\int_C\rho_{n}(z)\,z\,
    \pi(dz)\right)_+\right) ds\right]\\ 
&+&E\left[\frac{1}{h}\int_{t_n}^{(t_n+h)\wedge \theta_n}\psi\left(s,\rho_{n}(z)
     y_nY^{\rho_n,D_n}_{s^-}\right)-\psi\left(s, y_nY^{\rho_n,D_n}_{s^-}\right)\tilde
    \mu(ds,dz)\right]\geq \frac{\gamma_n}{h}.\\ 
\end{eqnarray*}
By the martingale's property, the third expectation on the left hand-side of the
last inequality vanishes and so we obtain 
\begin{eqnarray}\label{pm} 
& &E\left[\frac{1}{h}\int_{t_n}^{(t_n+h)\wedge \theta_n}\frac{\partial
    \psi}{\partial t}\left(s, y_nY^{\rho_n,D_n}_s\right) +
  A^{\rho_{n}}\left(s, y_nY^{\rho_n,D_n}_s,\psi, \Dy \psi\right)ds\right.\nonumber\\ 
&+&\left.
\frac{1}{h}\int_{t_n}^{(t_n+h)\wedge
    \theta_n} y_nY^{\rho_n,D_n}_s\left(\alpha-\beta+(\beta-\int_C\rho_{n}(z)\,z\,
    \pi(dz))_+\right) ds\right] 
\geq \frac{\gamma_n}{h}. 
\end{eqnarray}
From the definition of  $\gamma_n$, two cases are possible:\\
$\star$ Case $1$: if the set $\{n\geq 0: \gamma_n=0\}$ is finite, then there exists a
subsequence renamed $(\gamma_n)_{n\geq 0}$ such that $\gamma_n\neq 0$ for all
$n$ and we set $h=\sqrt{\gamma_n}$.\\
$\star$ Case $2$: if the set $\{n\geq 0: \gamma_n=0\}$ is not finite, then there
exists a subsequence renamed $(\gamma_n)_{n\geq 0}$ such that $\gamma_n= 0$ for all
$n$ and we set $h=n^{-1}$.\\ 
In both cases $\displaystyle{\frac{\gamma_n}{h}}\longrightarrow 0$ as $n$
tends to $\infty$ . We now send $n$ to infinity. The a.s. convergence of the
random value inside the expectation is obtained by the mean value
Theorem. Since $\int_C\pi(dz)<\infty$ and using the definition of $\theta_n$, the random variable
\begin{eqnarray*}
& &\frac{1}{h}\int_{t_n}^{(t_n+h)\wedge \theta_n}\frac{\partial
    \psi}{\partial t}\left(s, y_n Y^{\rho_n,D_n}_s\right) +
  A^{\rho_{n}}\left(s, y_nY^{\rho_n,D_n}_s,\psi, \Dy \psi\right)\\
&+& y_n Y^{\rho_n,D_n}_s\left(\alpha-\beta+(\beta-\int_C\rho_{ns}(z)\,z\,
    \pi(dz))_+\right) ds
\end{eqnarray*}
is essentially bounded, uniformly in $n$, on the stochastic interval 
$[t_n, (t_n+h)\wedge \theta_n]$. Sending $n$ to infinity, it follows by the dominated convergence theorem    
\begin{eqnarray*} 
\Dt \psi(t,y) + A^{\tilde \rho}(t,y,\psi, \Dy
\psi)+y\left(\alpha-\beta+(\beta-\int_C\tilde \rho(z)\,z\, \pi(dz))_+\right)\geq 0, 
\end{eqnarray*}
for all $\tilde \rho\in \Sigma$ and so
\begin{eqnarray}\label{pm0}
\Dt \psi(t,y) + H(t,y,\psi, \Dy \psi)\geq 0.
\end{eqnarray}
It remains to prove 
\begin{eqnarray*} 
-\Dy \psi(t,y) \geq 0. 
\end{eqnarray*}
For $\delta \in(0,1)$, we  set 
\begin{eqnarray*}
L_{s}=\left\{
\begin{array}{ll}
0 \mbox{ if } s=t_n^- \\
\delta \mbox{ if } s\geq t_n, 
\end{array}
\right.
\end{eqnarray*}
and so the process $D$ is a constant for all $s \geq t_n$.
We choose $\rho_s(z)=\tilde \rho(z)$ for all $s\geq t_n$ and $z\in C$, where $\tilde \rho\in \Sigma$.
We have $y_nY^{\rho,D}\in \Yc^0(t_n)$. Sending $n$ to infinity,
we have $y Y_{t^+}^{\rho,D}=y(1-\delta)$. Sending 
 $n\longrightarrow \infty$ in \reff{ppdgam}, by the dominated convergence theorem we get 
\begin{eqnarray*}
\psi(t,y)\leq \psi(t,y(1-\delta)).
\end{eqnarray*}
Sending $\delta \longrightarrow 0^+$, we obtain 
\begin{eqnarray}\label{sm}
-\Dy \psi(t,y) \geq 0.
\end{eqnarray}
Combining \reff{pm0} and (\ref{sm}), we conclude that $\tilde v$ is a
viscosity subsolution.\\
For supersolution inequality (\ref{superdef}), let $\psi \in C^1([0,T]\times
(0,\infty))$, $(\bar t,\bar y)\in [0,T)\times (0,\infty)$ such that $(\tilde v_{*}-\psi)(\bar t, \bar
y)=\min (\tilde v_{*}-\psi)=0$, we need to show
\begin{eqnarray}\label{superbar}
&&\min \left\{
\frac{\partial \psi(\bar t,\bar y)}{\partial t} + H_*\left(\bar t,\bar y,\psi, \Dy \psi\right),-
\Dy \psi(\bar t,\bar y) \right\} 
\leq 0. 
\end{eqnarray}
Suppose the contrary. Hence the left-hand side of \reff{superbar} is positive.
By smoothness of $\psi$ and since $H_*$ is lower semi-continuous, there exist
$\eta$ and $\epsilon$ satisfying:
\begin{eqnarray}\label{superbargamma}
& &\min \left\{
\frac{\partial \psi(t,y)}{\partial t} + H_*\left(t, y,\psi, \Dy
  \psi\right),-\Dy \psi(t,y) \right\}\geq \epsilon,
\end{eqnarray}
for all $(t,y)\in B_{\eta}(\bar t,\bar y)$, where $B_{\eta}(\bar t,\bar
y)=\{(t,y)\in [0,T]\times (0,\infty) \mbox{ such that }|t-\bar t|+|y-\bar y|\leq \eta \}$. By
changing $\eta$, we may assume that $B_{\eta}(\bar t,\bar y)\subset
[0,T)\times (0,\infty)$.\\
Since $(\bar t, \bar y)$ is a strict global minimizer of $\tilde v_{*} -\psi$
, there exists $\xi >0$ such that  
\begin{eqnarray*}
\Min_{(t,y)\in {\partial B_{\eta}}(\bar t,\bar y)}(\tilde v_{*} -\psi)(t,y)=\xi,
\end{eqnarray*} 
which implies $\tilde v_{*}(t,y) \geq \xi+\psi(t,y)$ for all $(t,y)\in
\partial B_\eta(\bar t,\bar y)$ the parabolic boundary of
$B_{\eta}(\bar t,\bar y)$.
From the definition of $\tilde v_{*}$, there exists a sequence $(t_n,y_n)\in
[0,T)\times (0,\infty) $  such that  $(t_n,y_n) \longrightarrow (\bar t,\bar y)$
and $\tilde v(t_n,y_n)\longrightarrow \tilde v_{*}(\bar t,\bar y)$ when
$n\longrightarrow \infty$. We suppose that $(t_n,y_n)\in B_{\eta}(\bar
t,\bar y) $.
Let $Y^{\rho,D}\in \Yc^0(t_n)$ be given and the stopping time $\theta_n$ defined by
\begin{eqnarray*}
\theta_n=\inf\{t\geq t_n\mbox{ such that }(t,y_nY_t^{\rho,D})\notin B_{\eta}(\bar
t,\bar y)\}\wedge T.
\end{eqnarray*}
Since the control $L\in \Lc_{t_n}$ is singular with possible jump at $t=t_n$, the couple
$(t,y_nY_t^{\rho,D})$ might jump out of $B_{\eta}(\bar t,\bar y)$ at $t_n$. If the control
$D$ makes alone the latter couple jump out of  $B_{\eta}(\bar t,\bar y)$, 
we set $\theta_n^D:=\theta_n$ else $\theta_n^D:=T$. 
In this case, the process $Y$ decreases. We know, from the dynamic programming principle
\reff{pppd} that 
\begin{eqnarray*}
\tilde v(t,y)\leq \tilde v(t,y(1-\delta))\,\mbox{ for all }0<\delta<1,
\end{eqnarray*}
and so the dual value
function $\tilde v$ in non-increasing with respect to $y$. However,
the point Poisson process could contributes to the jump out of
$B_{\eta}(\bar t,\bar y)$. In this case the dual value
function $\tilde v$ is not necessarily non-increasing in the direction of
the jump. The control $\rho$ could also contributes to hit the boundary of 
$B_{\eta}(\bar t,\bar y)$.  To overcome this problem, we introduce $\theta_j$ the first time
after $t_n$ the state process $Y$ jumps because of the point Poisson
process and  we set $\theta_n^\rho:=\theta_n$ when $y_nY^{\rho,D}$ jumps out 
$B_{\eta}(\bar t,\bar y)$ because of the control $\rho$ else $\theta_n^\rho:=T$. 
We set also $\theta_p:=\theta_n^\rho \wedge \theta_j$.  Note that, by right continuity of the paths, 
we have $\theta_p>t_n$ 
a.s. Let $\theta$ be the stopping time defined as follows $\theta:=\theta_n^D \wedge \theta_p$.\\
$\star$ On the set $\{\theta_n^D < \theta_p\}$. Let $(\theta_{n}^D,y^{'})$ be the intersection
between $\partial B_{\eta}(\bar t,\bar y)$ the parabolic boundary of
$B_{\eta}(\bar t,\bar y)$ and the line between
$(\theta_{n}^D,y_nY_{\theta_{n}^{D-}})$ and $(\theta_{n}^D,y_nY_{\theta_{n}^D})$ . From \reff{superbargamma}, 
we deduce that $\psi$ is non-increasing along this line in $\overline {B_{\eta}(\bar
  t,\bar y)}$. Since the dual value function $\tilde v$ is non-increasing
with respect to $y$, we have
\begin{eqnarray*}
\tilde v(\theta_{n}^D,y_nY_{\theta_{n}^D}) 
\geq \tilde v(\theta_{n}^D,y^{'}) 
\geq \psi(\theta_{n}^D,y^{'})+\xi
\geq \psi(\theta_{n}^D,y_nY_{\theta_{n}^{D-}})+\xi
\end{eqnarray*}
Using the inequality above and applying It\^o's formula to $\psi(t,y_nY_t^{\rho,D})$, we have
\begin{eqnarray}\label{itoy}
& & \tilde v(\theta_{n}^D,y_nY_{\theta_{n}^D}) \\
&\geq& \psi(\theta_{n}^D,y_nY_{\theta_{n}^{D-}})+\xi \nonumber\\
&\geq &\psi(t_n,y_n)+\int_{t_n}^{\theta_n^D}\Dt
  \psi(s,y_nY_s^{\rho,D}) + A^{\rho}(s,y_nY_s^{\rho,D},\psi, \Dy
  \psi)ds-\int_{t_n}^{\theta_n^D}\Dy \psi(s,y_nY_s^{\rho,D})y_nY^{\rho,D}_{s^-}dL_s\nonumber\\
&+& \int_{t_n}^{\theta_n^D}\int_C\psi\left(s,\rho_s(z)
   y_n Y^{\rho,D}_{s^-}\right)-\psi\left(s,y_nY^{\rho,D}_{s^-}\right)\tilde
    \mu(ds,dz)+\xi.\nonumber
\end{eqnarray}
For $t_n\leq s<\theta_n^D$, \reff{superbargamma} implies:
\begin{eqnarray}\label{contradiction1}
\Dt \psi(s,y_nY_s^{\rho,D})&+&
A^{\rho}(s, y_nY_s^{\rho,D},\psi, \Dy \psi)\nonumber\\
&+&y_nY_s^{\rho,D}\left(\alpha-\beta+(\beta-\int_C\rho_s(z)\,z\,
  \pi(dz))_+\right)\geq 0,
\end{eqnarray}
\begin{eqnarray}\label{contradiction2}
-\Dy \psi(s,y_nY_s^{\rho,D}) \geq 0.
\end{eqnarray}
Substituting \reff{contradiction1} and \reff{contradiction2} into \reff{itoy},
we have
\begin{eqnarray}\label{contradiction22}
& &\tilde v(\theta_{n}^D,y_nY_{\theta_{n}^D}^{\rho,D})+\int_{t_n}^{\theta_n^D}y_nY_s^{\rho,D}\left(\alpha-\beta+(\beta-\int_C\rho_s(z)\,z\,
  \pi(dz))_+\right)ds\\
&\geq& \psi(t_n,y_n)+\int_{t_n}^{\theta_n^D}\int_C\psi\left(s,\rho_s(z)
   y_n Y^{\rho,D}_{s^-}\right)-\psi\left(s,y_nY^{\rho,D}_{s^-}\right)\tilde
    \mu(ds,dz)
+\xi. \nonumber
\end{eqnarray} 
$\star$ On the set $\{\theta_n^D \geq \theta_p\}$, we have 
\begin{eqnarray}\label{itoy1}
 \tilde v(\theta_{p},y_nY_{\theta_{p}}) 
&\geq& \psi(\theta_{p},y_nY_{\theta_{p}}) \\
&\geq &\psi(t_n,y_n)+\int_{t_n}^{\theta_p}\Dt
  \psi(s,y_nY_s^{\rho,D}) + A^{\rho}(s,y_nY_s^{\rho,D},\psi, \Dy
  \psi)ds\nonumber \\
&-&\int_{t_n}^{\theta_p}\Dy \psi(s,y_nY_s^{\rho,D})y_nY^{\rho,D}_{s^-}dL_s\nonumber\\
&+& \int_{t_n}^{\theta_p}\int_C\psi\left(s,\rho_s(z)
   y_n Y^{\rho,D}_{s^-}\right)-\psi\left(s,y_nY^{\rho,D}_{s^-}\right)\tilde
    \mu(ds,dz).\nonumber
\end{eqnarray}
For $t_n\leq s<\theta_p$, \reff{superbargamma} implies:
\begin{eqnarray}\label{contradiction5}
\Dt \psi(s,y_nY_s^{\rho,D})&+&
A^{\rho}(s, y_nY_s^{\rho,D},\psi, \Dy \psi)\nonumber\\
&+&y_nY_s^{\rho,D}\left(\alpha-\beta+(\beta-\int_C\rho_s(z)\,z\,
  \pi(dz))_+\right)\geq \epsilon,
\end{eqnarray}
\begin{eqnarray}\label{contradiction6}
-\Dy \psi(s,y_nY_s^{\rho,D}) \geq 0.
\end{eqnarray}
Substituting \reff{contradiction5} and \reff{contradiction6} into \reff{itoy1},
we have
\begin{eqnarray}\label{contradiction23}
& &\tilde v(\theta_{p},y_nY_{\theta_{p}})+\int_{t_n}^{\theta_p}y_nY_s^{\rho,D}
\left(\alpha-\beta+(\beta-\int_C\rho_s(z)\,z\,\pi(dz))_+\right)ds\\
&\geq& \psi(t_n,y_n)+\int_{t_n}^{\theta_p}\int_C\Big(\psi\left(s,\rho_s(z)
    y_nY^{\rho,D}_{s^-}\right)-\psi\left(s,y_nY^{\rho,D}_{s^-}\right)\Big)\tilde
    \mu(ds,dz)+\epsilon (\theta_{p}-t_n). \nonumber
\end{eqnarray} 
Putting the two cases \reff{contradiction22} and \reff{contradiction23}
together, we get
\begin{eqnarray} \label{contradiction25}
& &E\left[\tilde
  v(\theta,y_nY_{\theta}^{\rho,D})+\int_{t_n}^{\theta}y_nY_s^{\rho,D}
\left(\alpha-\beta+(\beta-\int_C\rho_s(z)\,z\,\pi(dz))_+\right)ds\right]\\ 
&\geq&E\left[1_{\{\theta_n^D<\theta_p\} }\left(\tilde
  v(\theta_n^D,y_nY_{\theta_n^D}^{\rho,D})+\int_{t_n}^{\theta_n^D}
y_nY_s^{\rho,D}\left(\alpha-\beta+(\beta-\int_C\rho_s(z)\,z\,\pi(dz))_+\right)ds\right)
\right]\nonumber\\  
  &+&E\left[1_{\{\theta_n^D\geq\theta_p\} }\left(\tilde
  v(\theta_p,y_nY_{\theta_p}^{\rho,D})+\int_{t_n}^{\theta_p}y_nY_s^{\rho,D}
\left(\alpha-\beta+(\beta-\int_C\rho_s(z)\,z\,\pi(dz))_+\right)ds\right)\right]\nonumber\\
&\geq&\psi(t_n,y_n)+\xi P(\theta_n^D<\theta_p)+\epsilon 
E[1_{\{\theta_n^D\geq\theta_p\}}(\theta_p-t_n)]\nonumber
\end{eqnarray} 
Suppose that for all $\xi^{'}>0$, there exists $Y\in \Yc^0(t_n)$ such that 
\begin{eqnarray}\label{xip}
\xi P(\theta_n^D<\theta_p)+\epsilon 
E[1_{\{\theta_n^D\geq\theta_p\}}(\theta_p-t_n)] \leq \xi^{'}
\end{eqnarray} 
Since $\theta_p>t_n$ a.s. and $0\leq \epsilon E[1_{\{\theta_n^D\geq\theta_p\}}(\theta_p-t_n)]
\leq \xi^{'}$, for $\xi^{'}$ sufficiently small, we deduce that 
$\theta_n^D<\theta_p$ a.s. Inequality \reff{xip} implies $\xi\leq 0$ for
$\xi^{'}$ sufficiently small which is false and so there exists $\zeta>0$
such that for all $Y\in 
\Yc^0(t_n)$, we have
\begin{eqnarray}\label{xipj}
\xi P(\theta_n^D<\theta_p)+\epsilon 
E[1_{\{\theta_n^D\geq\theta_p\}}(\theta_p-t_n)] \geq \zeta.
\end{eqnarray}  
Inequalities \reff{contradiction25} and \reff{xipj} imply
\begin{eqnarray}\label{contradiction3}
& &\tilde v(t_n,y_n) +\delta_n\\
&\leq& E\left [
\tilde v(\theta,y_n Y_{\theta}^{\rho,D})
+\int_{t_n}^{\theta}y_nY_s^{\rho,D}\Big(\alpha-\beta+(\beta-\int_C\rho_s(z)\,z\,
\pi(dz))_+\Big)ds\right]-\zeta\nonumber
\end{eqnarray}
where $\delta_n:=\psi(t_n,y_n)-\tilde v(t_n,y_n) $.\\
Since $\delta_n=\psi(t_n,y_n)-\psi(\bar t,\bar y)+\tilde v_{*}(\bar t,\bar
y)-\tilde v(t_n,y_n)$, there exists $n_0\in \N$ such that for all $n\geq n_0
$, $\delta_n\geq -\frac{\zeta}{2}$. Inequality \reff {contradiction3} implies 
\begin{eqnarray*}
& &\tilde v(t_n,y_n) \nonumber\\
&\leq& \inf_{Y^{\rho,D}\in \Yc^0(t_n)}E\left[\tilde v(\theta,y_n Y_{\theta}^{\rho,D})
+\int_{t_n}^{\theta}y_nY_s^{\rho,D}\Big(\alpha-\beta+(\beta-\int_C\rho_s(z)\,z\,
\pi(dz))_+\Big)ds\right]- \frac{\zeta}{2}, 
\end{eqnarray*}
which is a contradiction with the dynamic programming principle and so we conclude 
that the dual value function  $\tilde v$ is a viscosity super-solution. 
\ep

\subsection{Proof of Theorem \ref{comparaison}}
For $\epsilon$, $\lambda$, $\delta$, $\zeta>0$, we define 
$\Phi :[0,T]\times (0,\infty)\times (0,\infty)\longrightarrow \R\cup \{-\infty\}$ as
\begin{eqnarray*}
\Phi(t,y_1,y_2)&:=&\tilde v_1^*(t,y_1)- \tilde
v_{2*}(t,y_2)-\frac{1}{\epsilon}(y_1-y_2)^2\\
&-& \delta \exp{(\lambda(T-t))} \left(y_1^{\gamma+1}+ y_2^{\gamma+1}\right)
- \zeta(\frac{1}{y_1^{\gamma+1}}+\frac{1}{y_2^{\gamma+1}}).
\end{eqnarray*}
Since $\tilde v_1^*$, $\tilde v_{2*}\in D_{\gamma}([0,T]\times (0,\infty))$, there
exists 
$(t^{*\epsilon,\delta,\lambda,\zeta},
  x^{*\epsilon,\delta,\lambda,\zeta},
  y^{*\epsilon,\delta,\lambda,\zeta})
\in   [0,T]\times (0,\infty)\times (0,\infty)$
which maximizes $\Phi$. By using the inequality
\begin{eqnarray*}
& &2\Phi(t^{*\epsilon,\delta,\lambda,\zeta},
y_1^{*\epsilon,\delta,\lambda,\zeta},
y_2^{*\epsilon,\delta,\lambda,\zeta})\\
&\geq&
\Phi(t^{*\epsilon,\delta,\lambda,\zeta},
y_1^{*\epsilon,\delta,\lambda,\zeta},
y_1^{*\epsilon,\delta,\lambda,\zeta})
+\Phi(t^{*\epsilon,\delta,\lambda,\zeta},
y_2^{*\epsilon,\delta,\lambda,\zeta},
y_2^{*\epsilon,\delta,\lambda,\zeta}), 
\end{eqnarray*}
we have 
\begin{eqnarray}\label{majxt}
\frac{2}{\epsilon}(y_1^{*\epsilon,\delta,\lambda,\zeta}-
y_2^{*\epsilon,\delta,\lambda,\zeta})^2
&\leq& 
\tilde v_1^*(t^{*\epsilon,\delta,\lambda,\zeta},
y_1^{*\epsilon,\delta,\lambda,\zeta})- 
\tilde v_1^*(t^{*\epsilon,\delta,\lambda,\zeta},
y_2^{*\epsilon,\delta,\lambda,\zeta})\nonumber\\
&+& \tilde v_{2*}(t^{*\epsilon,\delta,\lambda,\zeta},
y_1^{*\epsilon,\delta,\lambda,\zeta})- 
\tilde v_{2*}(t^{*\epsilon,\delta,\lambda,\zeta},
y_2^{*\epsilon,\delta,\lambda,\zeta}).\nonumber
\end{eqnarray}
since  $\tilde v^*_1$, $\tilde v_{2*}\in D_{\gamma}([0,T]\times (0,\infty))$, we have
\begin{eqnarray}\label{difference}
\frac{2}{\epsilon}
|y_1^{*\epsilon,\delta,\lambda,\zeta}
-y_2^{*\epsilon,\delta,\lambda,\zeta}|\leq
C(1+\frac{1}{(y_1^{*\epsilon,\delta,\lambda,\zeta})^{\gamma+1}}
+\frac{1}{(y_2^{*\epsilon,\delta,\lambda,\zeta})^{\gamma+1}}).
\end{eqnarray}
Using the inequality $\Phi(t^{*\epsilon,\delta,\lambda,\zeta},
y_1^{*\epsilon,\delta,\lambda,\zeta},
y_2^{*\epsilon,\delta,\lambda,\zeta})\geq \Phi(T,1,1)$ and 
since  $\tilde v^*_1$, $\tilde v_{2*} \in D_{\gamma}([0,T]\times (0,\infty))$, we have
\begin{eqnarray}\label{yz}
& &\delta \left((y_1^{*\epsilon,\delta,\lambda,\zeta})^{\gamma+1}
+(y_2^{*\epsilon,\delta,\lambda,\zeta})^{\gamma+1}\right)\\
&+&\zeta\left(\frac{1}{(y_1^{*\epsilon,\delta,\lambda,\zeta})^{*\gamma+1}}
+\frac{1}{(y_2^{*\epsilon,\delta,\lambda,\zeta})^{*\gamma+1}}\right)\nonumber\\
&\leq& C^{\delta,\zeta}\left(1+y_1^{*\epsilon,\delta,\lambda,\zeta}+
y_2^{*\epsilon,\delta,\lambda,\zeta}
+\frac{1}{(y_1^{*\epsilon,\delta,\lambda,\zeta})^\gamma}
+\frac{1}{(y_2^{*\epsilon,\delta,\lambda,\zeta})^\gamma}\right),\nonumber
\end{eqnarray}
where $C^{\delta,\zeta}$ is a constant depending only on $\delta$ and $\zeta$.
Inequality \reff{yz} implies either 
\begin{eqnarray}\label{separation}
\delta (y_1^{*\epsilon,\delta,\lambda,\zeta})^{\gamma+1}
+\zeta\frac{1}{(y_1^{*\epsilon,\delta,\lambda,\zeta})^{\gamma+1}}
\leq C^{\delta,\zeta}
\left(1+y_1^{*\epsilon,\delta,\lambda,\zeta}+\frac{1}{(y_1^{*\epsilon,\delta,\lambda,\zeta})^\gamma}\right) 
\end{eqnarray}
or 
\begin{eqnarray*}
\delta (y_2^{*\epsilon,\delta,\lambda,\zeta})^{\gamma+1}
+\tilde \zeta \frac{1}{(y_2^{*\epsilon,\delta,\lambda,\zeta})^{\gamma+1}}
\leq C^{\delta,\zeta}\left(1+y_2^{*\epsilon,\delta,\lambda,\zeta}
+\frac{1}{(y_2^{*\epsilon,\delta,\lambda,\zeta})^\gamma}\right). 
\end{eqnarray*}
Assume the first case, then there exist 
$M_1^{\delta,\zeta}$, 
$M_2^{\delta,\zeta}>0$ depending only on  $ \delta$ and $\zeta$ such
that  
$M_1^{ \delta,\zeta}
\leq y_1^{*\epsilon,\delta,\lambda,\zeta}
\leq M_2^{\delta,\zeta}$. 
Using Inequality \reff{yz}, we obtain
\begin{eqnarray*}
 \delta 
(y_2^{*\epsilon,\delta,\lambda,\zeta})^{\gamma+1}
\leq
C^{\delta,\zeta}
(1+y_1^{*\epsilon,\delta,\lambda,\zeta}
+y_2^{*\epsilon,\delta,\lambda,\zeta}
+\frac{1}{(y_1^{*\epsilon,\delta,\lambda,\zeta})^\gamma}
+\frac{1}
{(y_2^{*\epsilon,\delta,\lambda,\zeta})^\gamma})
\end{eqnarray*}
and
\begin{eqnarray*}
\zeta\frac{1}{(y_2^{*\epsilon,\delta,\lambda,\zeta})^{\gamma+1}} 
\leq C^{\delta,\zeta}
(1+y_1^{*\epsilon,\delta,\lambda,\zeta}
+y_2^{*\epsilon,\delta,\lambda,\zeta}
+\frac{1}{(y_1^{*\epsilon,\delta,\lambda,\zeta})^\gamma}
+\frac{1}
{(y_2^{*\epsilon,\delta,\lambda,\zeta})^\gamma}),
\end{eqnarray*}
which implies that  
$M_1^{\delta, \zeta}
\leq y_2^{*\epsilon,\delta,\lambda,\zeta}
\leq M_2^{\delta,\zeta}$. 
Since $y_1^{*\epsilon,\delta,\lambda,\zeta}$ and $y_2^{*\epsilon,\delta,\lambda,\zeta}$ 
are bounded from below, inequality \reff{difference} implies
\begin{eqnarray}\label{diff}
\big|y_1^{*\epsilon,\delta,\lambda,\zeta}-y_2^{*\epsilon,\delta,\lambda,\zeta}\big|\leq C_1\epsilon,
\end{eqnarray}
where $C_1$ is a positive constant independent of $\epsilon$. 
Using the boundedness of $y_1^{*\epsilon,\delta,\lambda,\zeta}$ 
and $y_2^{*\epsilon,\delta,\lambda,\zeta}$ and  \reff{diff},
along a subsequence $(t^{*\epsilon,\delta,\lambda,\zeta},
y_1^{*\epsilon,\delta,\lambda,\zeta},
y_2^{*\epsilon,\delta,\lambda,\zeta})$ 
converges when $\epsilon \longrightarrow 0$. Let's denote 
$(t^{*\delta,\lambda,\zeta},
y^{*\delta,\lambda,\zeta},
y^{*\delta,\lambda,\zeta})$  its limit.\\
From the definition of $t^{*\epsilon,\delta,\lambda,\zeta}$, two cases are possible:\\
$\star $ Case $1$: If the set $\{\epsilon
>0:\,t^{*\epsilon,\delta,\lambda,\zeta}=T \}$ is not finite, then there exists a
subsequence renamed $(t^{*\epsilon,\delta,\lambda,\zeta})_{\epsilon}$ such that 
$t^{*\epsilon,\delta,\lambda,\zeta}=T$. From inequality
$$\Phi(t,y,y)\leq
\Phi(t^{*\epsilon,\delta,\lambda,\zeta},
y_1^{*\epsilon,\delta,\lambda,\zeta},
y_2^{*\epsilon,\delta,\lambda,\zeta})$$ 
and since $\Phi$ is upper semi-continuous, we deduce
that 
\begin{eqnarray*}
\Phi(t,y,y)&\leq& 
\Limsup_{\epsilon \rightarrow 0}
\Phi(t^{*\epsilon,\delta,\lambda,\zeta},
y_1^{*\epsilon,\delta,\lambda,\zeta},
y_2^{*\epsilon,\delta,\lambda,\zeta})\\
&\leq&\Phi(T,y^{*\delta,\lambda,\zeta},y^{*\delta,\lambda,\zeta}),
\end{eqnarray*}
which implies
\begin{eqnarray*}
& &\tilde v_1^*(t,y)- \tilde v_{2*}(t,y)
-2 \delta \exp{(\lambda(T-t))}y^{\gamma+1}-2\zeta \frac{1}{y^{\gamma+1}}\\ 
&\leq& \tilde v_1^*(T,y^{*\delta,\lambda,\zeta})
- \tilde v_{2*}(T,y^{*\delta,\lambda,\zeta}).
\end{eqnarray*}
Using inequality $ \tilde v_1^*(T, y^{*\delta,\lambda,\zeta})
\leq \tilde v_{2*}(T,y^{*\delta,\lambda,\zeta})$ and sending $\lambda$, $\delta$, $\zeta$ $\longrightarrow 0^+$, 
we have
\begin{eqnarray*}
\tilde u^*(t,y) \leq \tilde v_*(t,y),\,\mbox{ for all }(t,y)\in [0,T]\times (0,\infty).
\end{eqnarray*}
$\star$ Case $2$: If the set $\{\epsilon >0:\,t^{*\epsilon,\delta,\lambda,\zeta}=T \}$ is finite, then there exists a
subsequence renamed $(t^{*\epsilon,\delta,\lambda,\zeta})_{\epsilon}$ such that 
$t^{*\epsilon,\delta,\lambda,\zeta}<T$. Our aim is to construct a regular function denoted $\tilde \psi_1$ 
(resp. $\tilde \psi_2$) satisfying  inequality \reff{superdef} (resp. \reff{subdef}).
We define $\psi_1$ and $\psi_2$ as follows 
\begin{eqnarray*}
\psi_1(t,y)&:=&
\tilde v_{*2}(t^{*\epsilon,\delta,\lambda,\zeta},y_2^{*\epsilon,\delta,\lambda,\zeta})
+\frac{1}{\epsilon}(y-y_2^{*\epsilon,\delta,\lambda,\zeta})^2\\
&+&\delta\left( \exp{( \lambda(T-t))}y^{\gamma+1}
+\exp{(\lambda(T-t^{*\epsilon,\delta,\lambda,\zeta}))}
(y_2^{*\epsilon,\delta,\lambda,\zeta})^{\gamma+1}\right)\\
&+&\zeta(\frac{1}{y^{\gamma+1}}+\frac{1}{(y_2^{*\epsilon,\delta,\lambda,\zeta})
^{\gamma+1}})+\Phi(t^{*\epsilon,\delta,\lambda,\zeta},
y_1^{*\epsilon,\delta,\lambda,\zeta},y_2^{*\epsilon,\delta,\lambda,\zeta}),\,\,(t,y)\in [0,T]\times (0,\infty),
\end{eqnarray*}
and
\begin{eqnarray*}
\psi_2(t,y)&:=&\tilde v^*_1(t^{*\epsilon,\delta,\lambda,\zeta},
y_1^{*\epsilon,\delta,\lambda,\zeta})
-\frac{1}{\epsilon}(y_1^{*\epsilon,\delta,\lambda,\zeta}-y)^2\\
&-&\delta\left( \exp{(\lambda(T-t^{*\epsilon,\delta,\lambda,\zeta}))}
(y_1^{*\epsilon,\delta,\lambda,\zeta})^{\gamma+1}+\exp{( \lambda(T-t))}y^{\gamma+1}\right)\\
&-& \zeta(\frac{1}{(y_1^{*\epsilon,\delta,\lambda,\zeta})^{\gamma+1}}+\frac{1}{y^{\gamma+1}})
-\Phi(t^{*\epsilon,\delta,\lambda,\zeta},
y_1^{*\epsilon,\delta,\lambda,\zeta},y_2^{*\epsilon,\delta,\lambda,\zeta})
,\,\,(t,y)\in [0,T]\times (0,\infty).
\end{eqnarray*}
From inequalities \reff{minorant4} and \reff{minorant5}, we have
\begin{eqnarray}\label{bornetildev}
\tilde U(y)-|K^{'}| y \leq \tilde v(t,y)\leq \tilde U(y)+|K| y.
\end{eqnarray}
We define $\tilde \psi_1$ and $\tilde \psi_2$ as follows
\begin{eqnarray}\label{sub-1}
\tilde \psi_1(t,y)=\left \{ \begin{array}{ll}
2\tilde U(y)&\,\mbox{ for all }0<y\leq \underline  y_1\wedge \frac{M_1^{\delta,\zeta}}{2} ,\\
\psi_1(t,y)&\,\mbox{ for all } M_1^{\delta,\zeta}\leq y\leq M_2^{\delta,\zeta}\\
M_1 y&\,\mbox{ for all }y\geq  2 M_2^{\delta,\zeta},
\end{array}
\right.
\end{eqnarray}
and
\begin{eqnarray}\label{sub-}
\tilde \psi_2(t,y)=\left \{ \begin{array}{ll}
\frac{\tilde U(y)}{2}&\,\mbox{ for all }0<y\leq \underline  y_2\wedge \frac{M_1^{\delta,\zeta}}{2} ,\\
\psi_2(t,y)&\,\mbox{ for all } M_1^{\delta,\zeta}\leq y\leq M_2^{\delta,\zeta}\\
M_2 y+M_2^{'}&\,\mbox{ for all }y\geq  2 M_2^{\delta,\zeta},
\end{array}
\right.
\end{eqnarray}
where $\underline  y_1= \big(\gamma|K|\big)^{-\frac{1}{\gamma+1}}$, $M_1$ satisfies
\begin{eqnarray}\label{croi1}
M_1\geq \frac{ (2 M_2^{\delta,\zeta})^{-(\gamma +1)}}{\gamma}+|K|
\vee \frac {\partial \psi_1 (t^{*\epsilon,\delta,\lambda,\zeta},y_1^{*\epsilon,\delta,\lambda,\zeta} )}{\partial y} ,
\end{eqnarray}
$\underline  y_2= \big(2\gamma|K^{'}|\big)^{-\frac{1}{\gamma+1}}$, $M_2$  satisfies
\begin{eqnarray}\label{choixm2}
\frac {\partial \psi_2 (t^{*\epsilon,\delta,\lambda,\zeta},y_2^{*\epsilon,\delta,\lambda,\zeta} )}{\partial y}
\leq M_2\leq \tilde v_{2*d}^{'}(t^{*\epsilon,\delta,\lambda,\zeta},2 M_2^{\delta,\zeta} ),
\end{eqnarray} 
where $\tilde v_{2*d}^{'}$  is the right-hand derivative of $\tilde v_{2*}$ with respect to the variable $y$ and 
$M_2^{'}=\tilde v_{*2}(t,2 M_2^{\delta,\zeta})- 2 M_2  M_2^{\delta,\zeta}-1$.
The choice of $M_2$ is possible since  $\tilde v_{2*}$ is convex and nonincreasing 
( see assumptions of the comparison theorem ) and so 
\begin{eqnarray*}
\tilde v_{2*d}^{'}(t^{*\epsilon,\delta,\lambda,\zeta},y_2^{*\epsilon,\delta,\lambda,\zeta}  ) 
\leq \tilde v_{2*d}^{'}(t^{*\epsilon,\delta,\lambda,\zeta},2 M_2^{\delta,\zeta} )
\end{eqnarray*}
and from the optimality of 
$(t^{*\epsilon,\delta,\lambda,\zeta},y_2^{*\epsilon,\delta,\lambda,\zeta} )$, we have
\begin{eqnarray}\label{croi2}
\frac {\partial \psi_2 (t^{*\epsilon,\delta,\lambda,\zeta},y_2^{*\epsilon,\delta,\lambda,\zeta} )}{\partial y}\leq \tilde v_{2*d}^{'}(t^{*\epsilon,\delta,\lambda,\zeta},y_2^{*\epsilon,\delta,\lambda,\zeta}  ).
\end{eqnarray}
From  the definition of $\underline  y_1$ and inequalities \reff{bornetildev}, \reff{croi1}, we have
\begin{eqnarray*}
\left \{ \begin{array}{ll}
\tilde v(t,y)\leq 2\tilde U(y) &\mbox{ for all }0<y\leq \underline  y_1 \\
\tilde v(t,y)\leq \tilde U(y)+|K|y\leq M_1y&\mbox{ for all }y\geq  2 M_2^{\delta,\zeta}
\end{array}
\right.
\end{eqnarray*}
and so one could obtain  inequality \reff{superdef} for $\tilde \psi_1$.\\ 
From  the definition of $\underline  y_2$, inequalities \reff{bornetildev}, \reff{croi2} and using the convexity of 
$ \tilde v_{2*d}$, we have
\begin{eqnarray*}
\left \{ \begin{array}{ll}
\tilde v(t,y)\geq \frac{\tilde U(y)}{2} &\mbox{ for all }0<y\leq \underline  y_2 \\
\tilde v(t,y)\geq M_2y+M_2^{'}&\mbox{ for all }y\geq  2 M_2^{\delta,\zeta}
\end{array}
\right.
\end{eqnarray*}
and so one could obtain  inequality \reff{subdef} for $\tilde \psi_2$.\\
To prove the comparison theorem, we need to derive an equivalent formulation of viscosity solutions 
an in Soner \cite{son86b} Lemma 2.1. For this, we show that the control $\rho$ runs along a compact 
set and the lower semi-continuous envelope of $H$ is continuous in its arguments which is the object 
of the next lemma. We denote by $\bar \rho:=\frac{\beta}{\Min_{1\leq i\leq d}\delta_i\pi_i}$. 
We define $\Sigma^{'}$ by 
\begin{eqnarray}\label{sigmap}
\Sigma^{'}=\left\{ \rho =(\rho_i)_{1\leq i\leq d},\,\, 
0\leq\rho_i\leq\frac{ 2 M_2^{\delta,\zeta}+1}{y_1^{*\epsilon,\delta,\lambda,\zeta}}
\vee \frac{ 2 M_2^{\delta,\zeta}+1}{y_2^{*\epsilon,\delta,\lambda,\zeta}}
\vee  \bar \rho \right\}
\end{eqnarray}
and the Hamiltonian $\tilde H$ by
\begin{eqnarray}\label{tildeH}
\tilde H\left(t,y,\tilde v, \Dy {\tilde v}\right)
=\Inf_{\rho \in \Sigma^{'}}\left\{A^\rho\left(t,y,\tilde v, \Dy {\tilde v}\right)+
y\left(\alpha-\beta+(\beta-\Sum_{i=1}^d\rho_i \delta_i \pi_i)_+\right)\right\},
\end{eqnarray} 
\begin{Lemma}
We assume that $\tilde v_{2*}$ is convex and nonincreasing. Then,
we have the following inequalities
\begin{eqnarray}\label{comparaison1}
& &\min \left\{
\frac{\partial \psi_1}{\partial t}
(t^{*\epsilon,\delta,\lambda,\zeta},
y_1^{*\epsilon,\delta,\lambda,\zeta}) 
+\tilde H\left(t^{*\epsilon,\delta,\lambda,\zeta},
y_1^{*\epsilon,\delta,\lambda,\zeta},
\tilde v_1^*, 
\Dy {\psi_1}\right),\right.\\
& &\left.-\Dy {\psi_1}(t^{*\epsilon,\delta,\lambda,\zeta},
y_1^{*\epsilon,\delta,\lambda,\zeta}) \right\}\geq 0,\nonumber 
\end{eqnarray} 
and
\begin{eqnarray}\label{comparaison2}
&&\min \left\{
\frac{\partial \psi_2}{\partial
  t}(t^{*\epsilon,\delta,\lambda,\zeta},y_2^{*\epsilon,\delta,\lambda,\zeta}) 
+\tilde H \left(t^{*\epsilon,\delta,\lambda,\zeta},
y_2^{*\epsilon,\delta,\lambda,\zeta},\tilde v_{2*}, \Dy {\psi_2}\right),\right.\\
&&\left.-\Dy {\psi_2}(t^{*\epsilon,\delta,\lambda,\zeta},
y_2^{*\epsilon,\delta,\lambda,\zeta}) \right\}\leq 0.\nonumber 
\end{eqnarray}  
\end{Lemma}
{\bf Proof.}
If there exists $0\leq i_0\leq d$ such that $\displaystyle{\rho_{i_0}\geq \bar \rho}$,
then $(\beta-\Sum_{i=1}^d\rho_i \delta_i \pi_i)_+=0$. 
Let $\rho\in \Sigma$ be a fixed vector. If $\rho_i \geq \frac{ 2 M_2^{\delta,\zeta}}{y_1^{*\epsilon,\delta,\lambda,\zeta}} 
\vee  \bar \rho$, $0\leq i \leq d$ , then we have
\begin{eqnarray*}
F_1(\rho_i)&:=&\tilde \psi_1(t^{*\epsilon,\delta,\lambda,\zeta},\rho_i y_1^{*\epsilon,\delta,\lambda,\zeta})- 
\tilde \psi_1 (t^{*\epsilon,\delta,\lambda,\zeta},y_1^{*\epsilon,\delta,\lambda,\zeta})\\
&-&(\rho_i -1) y_1^{*\epsilon,\delta,\lambda,\zeta} 
\Dy {\tilde \psi_1}(t^{*\epsilon,\delta,\lambda,\zeta},y_1^{*\epsilon,\delta,\lambda,\zeta})\\
&=&(M_1- \Dy {\psi_1}(t^{*\epsilon,\delta,\lambda,\zeta},y_1^{*\epsilon,\delta,\lambda,\zeta}))
\rho_i y_1^{*\epsilon,\delta,\lambda,\zeta}
-\psi_1 (t^{*\epsilon,\delta,\lambda,\zeta},y_1^{*\epsilon,\delta,\lambda,\zeta})\\
&+&y_1^{*\epsilon,\delta,\lambda,\zeta} \Dy {\psi_1}(t^{*\epsilon,\delta,\lambda,\zeta},y_1^{*\epsilon,\delta,\lambda,\zeta}).
\end{eqnarray*}
From \reff{croi1}, we deduce that the function $\rho_i \longrightarrow F_1(\rho_i)$ is non-decreasing and so
\begin{eqnarray}\label{ham1}
& &\Inf_{\rho \in \Sigma}
\left\{A^\rho\left(t^{*\epsilon,\delta,\lambda,\zeta},
y_1^{*\epsilon,\delta,\lambda,\zeta},\tilde \psi_1, \Dy {\tilde \psi_1}\right)+
y_1^{*\epsilon,\delta,\lambda,\zeta}\left(\alpha-\beta+(\beta-\Sum_{i=1}^d\rho_i \delta_i \pi_i)_+\right)\right\}\nonumber\\
&=&\Inf_{\rho \in \Sigma^{'}}
\left\{A^\rho\left(t^{*\epsilon,\delta,\lambda,\zeta},y_1^{*\epsilon,\delta,\lambda,\zeta},\tilde \psi_1, 
\Dy {\tilde \psi_1}\right)+
y_1^{*\epsilon,\delta,\lambda,\zeta}\left(\alpha-\beta+(\beta-\Sum_{i=1}^d\rho_i \delta_i \pi_i)_
+\right)\right\}\nonumber\\
&:=&\Inf_{\rho \in \Sigma^{'}}f_1(t^{*\epsilon,\delta,\lambda,\zeta},y_1^{*\epsilon,\delta,\lambda,\zeta},\rho )\nonumber\\
&:=&v_1^{opt}(t^{*\epsilon,\delta,\lambda,\zeta},y_1^{*\epsilon,\delta,\lambda,\zeta}).
\end{eqnarray}
The criterion of the optimization problem \reff{ham1} is continuous with respect to $\rho$, $\Sigma^{'}$ 
is a compact and so there exists  $\rho^*_1$ solution of \reff{ham1}. 
We consider a sequence $(t_k,y_k)_{k}\in [0,T]\times (0,\infty)$ such that $(t_k,y_k)\longrightarrow(t^{*\epsilon,\delta,\lambda,\zeta},y_1^{*\epsilon,\delta,\lambda,\zeta}) $ when $n$ goes to infinity. We denote by $\rho_k$ the optimum i.e. 
\begin{eqnarray}\label{conH1}
v_1^{opt}(t_k,y_k)=f_1(t_k,y_k,\rho_k ).
\end{eqnarray}
Since $\rho_k\in \Sigma^{'} $ which is compact, then along a subsequence denoted also by $(\rho_k)_k$, we have
$\rho_k\longrightarrow \bar \rho$. From the Taylor expansion formula and using the continuity of $f_1$ in her arguments, we have
\begin{eqnarray}\label{conH2}
f_1(t_k,y_k,\rho_k )&=& f_1(t^{*\epsilon,\delta,\lambda,\zeta},y_1^{*\epsilon,\delta,\lambda,\zeta},\bar \rho )+o(1)\nonumber\\
&\geq&v_1^{opt}(t^{*\epsilon,\delta,\lambda,\zeta},y_1^{*\epsilon,\delta,\lambda,\zeta})+o(1)
\end{eqnarray}
From \reff{conH1} and \reff{conH2}, we deduce that
\begin{eqnarray}
v_1^{opt}(t_k,y_k)&\geq&v_1^{opt}(t^{*\epsilon,\delta,\lambda,\zeta},y_1^{*\epsilon,\delta,\lambda,\zeta})+o(1).
\end{eqnarray}
To obtain the converse inequality, we have
\begin{eqnarray*}
v_1^{opt}(t_k,y_k)&\leq& f_1(t_k,y_k, \rho^*_1  )\\
&=&f_1(t^{*\epsilon,\delta,\lambda,\zeta},y_1^{*\epsilon,\delta,\lambda,\zeta}, \rho^*_1  )+o(1)\\
&=&v_1^{opt}(t^{*\epsilon,\delta,\lambda,\zeta},y_1^{*\epsilon,\delta,\lambda,\zeta})+o(1),
\end{eqnarray*}
and so $\Lim_{k\longrightarrow \infty}v_1^{opt}(t_k,y_k)=v_1^{opt}(t^{*\epsilon,\delta,\lambda,\zeta},y_1^{*\epsilon,\delta,\lambda,\zeta})$. This proves that $v_1^{opt}$ is continuous in $
(t^{*\epsilon,\delta,\lambda,\zeta},y_1^{*\epsilon,\delta,\lambda,\zeta})$ and so
\begin{eqnarray}\label{H*}
H_*\left(t^{*\epsilon,\delta,\lambda,\zeta},y_1^{*\epsilon,\delta,\lambda,\zeta},\tilde \psi_1, \Dy {\tilde \psi_1} \right)
=\tilde H\left(t^{*\epsilon,\delta,\lambda,\zeta},y_1^{*\epsilon,\delta,\lambda,\zeta},\tilde \psi_1, \Dy {\tilde \psi_1} \right).
\end{eqnarray} 
From equality \reff{sub-1}, the function $\tilde v^*_1- \tilde \psi_1$ has a strict global minimum at 
$(t^{*\epsilon,\delta,\lambda,\zeta},y_1^{*\epsilon,\delta,\lambda,\zeta})\in [0,T)\times (0,\infty)$. 
Using the definition of viscosity supersolutions (see inequality \reff{superdef}), equation \reff{H*} and \reff{ham1}, we obtain
\begin{eqnarray}\label{sub1}
&&\min \Big\{
\frac{\partial \tilde \psi_1}{\partial t}(t^{*\epsilon,\delta,\lambda,\zeta},y_1^{*\epsilon,\delta,\lambda,\zeta}) +
\tilde H
\Big(t^{*\epsilon,\delta,\lambda,\zeta},
y_1^{*\epsilon,\delta,\lambda,\zeta},\tilde \psi_1, \frac{\partial \tilde \psi_1}{\partial y} \Big),\nonumber \\
&-&\frac{\partial \tilde \psi_1}{\partial y}(t^{*\epsilon,\delta,\lambda,\zeta},y_1^{*\epsilon,\delta,\lambda,\zeta}) \Big\}
\geq 0, 
\end{eqnarray} 
Let $\rho\in \Sigma$ 
be a fixed vector. If $\rho_i \geq \frac{ 2 M_2^{\delta,\zeta}}{y_2^{*\epsilon,\delta,\lambda,\zeta}} 
\vee  \bar \rho$, then we have
\begin{eqnarray*}
& &F_2(\rho_i)\\
&: =&\tilde \psi_2(t^{*\epsilon,\delta,\lambda,\zeta},\rho_i y_2^{*\epsilon,\delta,\lambda,\zeta})- 
\tilde \psi_2 (t^{*\epsilon,\delta,\lambda,\zeta},y_2^{*\epsilon,\delta,\lambda,\zeta})
-(\rho_i -1) y_2^{*\epsilon,\delta,\lambda,\zeta} 
\Dy {\tilde \psi_2}(t^{*\epsilon,\delta,\lambda,\zeta},y_2^{*\epsilon,\delta,\lambda,\zeta})\\
&=&(M_2- \Dy {\psi_2}(t^{*\epsilon,\delta,\lambda,\zeta},y_2^{*\epsilon,\delta,\lambda,\zeta}))
\rho_i y_2^{*\epsilon,\delta,\lambda,\zeta}
-\psi_2 (t^{*\epsilon,\delta,\lambda,\zeta},y_2^{*\epsilon,\delta,\lambda,\zeta})\\
&+&y_2^{*\epsilon,\delta,\lambda,\zeta} \Dy {\psi_2}(t^{*\epsilon,\delta,\lambda,\zeta},y_2^{*\epsilon,\delta,\lambda,\zeta}).
\end{eqnarray*}
From the definition of $\tilde \psi_2$ and inequality \reff{choixm2}, we deduce that the function $\rho_i \longrightarrow F_2(\rho_i)$ is non-decreasing and so
\begin{eqnarray}\label{ham}
& &\Inf_{\rho \in \Sigma}
\left\{A^\rho\left(t^{*\epsilon,\delta,\lambda,\zeta},
y_2^{*\epsilon,\delta,\lambda,\zeta},\tilde \psi_2, \Dy {\tilde \psi_2}\right)+
y_2^{*\epsilon,\delta,\lambda,\zeta}\left(\alpha-\beta+(\beta-\Sum_{i=1}^d\rho_i \delta_i \pi_i)_+\right)\right\}\\
&=&\Inf_{\rho \in \Sigma^{'}}
\left\{A^\rho\left(t^{*\epsilon,\delta,\lambda,\zeta},y_2^{*\epsilon,\delta,\lambda,\zeta}
,\tilde \psi_2, \Dy {\tilde \psi_2}\right)+
y_2^{*\epsilon,\delta,\lambda,\zeta}\left(\alpha-\beta+(\beta-\Sum_{i=1}^d\rho_i \delta_i \pi_i)_+\right)\right\},\nonumber
\end{eqnarray}
From equality \reff{sub-},$\tilde v_{2*}-
 \tilde \psi_2$ has a  strict global maximum at 
$(t^{*\epsilon,\delta,\lambda,\zeta},y_2^{*\epsilon,\delta,\lambda,\zeta})\in [0,T)\times (0,\infty)$. 
Using the definition of viscosity sub-solutions (see inequality \reff{subdef}) and \reff{ham}, we obtain
\begin{eqnarray}\label{sub2}
&&\min \Big\{
\frac{\partial \tilde \psi_2}{\partial t}(t^{*\epsilon,\delta,\lambda,\zeta},y_2^{*\epsilon,\delta,\lambda,\zeta})
+
\tilde H
\Big(t^{*\epsilon,\delta,\lambda,\zeta},
y_2^{*\epsilon,\delta,\lambda,\zeta},\tilde \psi_2, \frac{\partial \tilde \psi_2}{\partial y} \Big), \nonumber\\
&&-\frac{\partial \tilde \psi_2}{\partial y}(t^{*\epsilon,\delta,\lambda,\zeta},y_2^{*\epsilon,\delta,\lambda,\zeta}) \Big\} \geq 0. 
\end{eqnarray} 
From inequalities \reff{sub1} and \reff{sub2}, using the fact the control $\rho$ runs along a compact set and 
repeating arguments of Soner \cite{son86b} Lemma 2.1, we easily obtain an equivalent formulation of viscosity solutions given by inequalities \reff{comparaison1} and \reff{comparaison2}.
\ep\\
\\
We come back to the proof of the comparison theorem.
Remarking that $\min\{d,e\}-\min\{f,g\}\geq 0$ implies either $d-f\geq 0$ 
or $e-g\geq 0$, inequalities \reff{comparaison1} and \reff{comparaison2} imply 
\begin{eqnarray}\label{uniqueness1} 
- \delta \lambda
\exp{(\lambda(T-t^{*\epsilon,\delta,\lambda,\zeta}))}
\left(
(y_1^{*\epsilon,\delta,\lambda,\zeta})^{\gamma+1}+
(y_2^{*\epsilon,\delta,\lambda,\zeta})^{\gamma+1}\right)
+T_1-T_2\geq 0,
\end{eqnarray}
or 
\begin{eqnarray}\label{uniqueness2}
& &-\delta \exp{\left( \lambda(T-t^{*\epsilon,\delta,\lambda,\zeta})\right)}
\left((y_1^{*\epsilon,\delta,\lambda,\zeta})^\gamma
+(y_2^{*\epsilon,\delta,\lambda,\zeta})^{\gamma}\right)\\
&+& \zeta\left(\frac{1}{(y_1^{*\epsilon,\delta,\lambda,\zeta})^{\gamma+2}}
+\frac{1}{(y_2^{*\epsilon,\delta,\lambda,\zeta})^{\gamma+2}}\right)\geq 0,\nonumber
\end{eqnarray}
where
\begin{eqnarray}\label{criter}
&T_1&:=\Inf_{\rho \in \Sigma^{'}}
\left\{
\Sum_{i=1}^d 
\pi_i
\Big(\tilde v_1^*(t^{*\epsilon,\delta,\lambda,\zeta},\rho_i 
y_1^{*\epsilon,\delta,\lambda,\zeta})
-\tilde v_1^*(t^{*\epsilon,\delta,\lambda,\zeta},
y_1^{*\epsilon,\delta,\lambda,\zeta})
\right.
\\
&-&(\rho_i-1)y_1^{*\epsilon,\delta,\lambda,\zeta}
\Big(\frac{2}{\epsilon}
(y_1^{*\epsilon,\delta,\lambda,\zeta}
-y_2^{*\epsilon,\delta,\lambda,\zeta})+ \delta(\gamma+1)
\exp{( \lambda(T-t^{*\epsilon,\delta,\lambda,\zeta}))}
(y_1^{*\epsilon,\delta,\lambda,\zeta})^\gamma\nonumber\\
&-&
\left.  \zeta(\gamma+1)\frac{1}{(y_1^{*\epsilon,\delta,\lambda,\zeta})^{\gamma
  +2}}
\Big)
\Big)
+y_1^{*\epsilon,\delta,\lambda,\zeta}
\left(
\alpha-\beta+(\beta-\Sum_{i=1}^{d}\rho_i\delta_i\pi_i)_+
\right)
\right\},\nonumber
\end{eqnarray}
\begin{eqnarray}
&T_2&:=\Inf_{\rho \in  \Sigma^{'}}\left\{\Sum_{i=1}^d 
\pi_i
\Big(
\tilde v_{2*}(t^{*\epsilon,\delta,\lambda,\zeta},\rho_i
y_2^{*\epsilon,\delta,\lambda,\zeta})
-\tilde v_{2*}(t^{*\epsilon,\delta,\lambda,\zeta},
y_2^{*\epsilon,\delta,\lambda,\zeta})\right.\nonumber\\
&-&(\rho_i-1)y_2^{*\epsilon,\delta,\lambda,\zeta}
\Big(
\frac{2}{\epsilon}(y_1^{*\epsilon,\delta,\lambda,\zeta}
-y_2^{*\epsilon,\delta,\lambda,\zeta})- \delta(\gamma+1)
\exp{(\lambda(T-t^{*\epsilon,\delta,\lambda,\zeta}))}
(y_2^{*\epsilon,\delta,\lambda,\zeta})^\gamma\nonumber\\
&+& \left. \zeta(\gamma+1)\frac{1}
{(y_2^{*\epsilon,\delta,\lambda,\zeta})^{\gamma+2}}
\Big)
\Big)
+y_2^{*\epsilon,\delta,\lambda,\zeta}\left(\alpha-\beta
+(\beta-\Sum_{i=1}^{d}\rho_i\delta_i\pi_i)_+\right)\right\}\nonumber
\end{eqnarray}
and $\Sigma^{'}$ is defined in \reff{sigmap}. From Lemma \ref{conv}
$\tilde  v_{2*}$ is continuous w.r.t the state variable ($\tilde  v_{2*}$ is convex on $(0,\infty)$) and so
the criterion of the optimization problem \reff{criter} is continuous  w.r.t $\rho$. The set $\Sigma^{'}$ is compact 
and since $\Lim_{y\longrightarrow 0} v_{2*}(t,y)=\infty$, 
there exists a solution denoted by $\rho^{*\epsilon,\delta,\lambda,\zeta}$, to the optimization problem \reff{criter} 
satisfying $\rho^{*\epsilon,\delta,\lambda,\zeta}_i>0$, for all $1\leq i\leq d$. \\
We define $f$ as follows:
\begin{eqnarray*}
f(\rho)&:=&\Sum_{i=1}^d 
\pi_i
\Big(\tilde v_{2*}(t^{*\epsilon,\delta,\lambda,\zeta},\rho_i 
y_2^{*\epsilon,\delta,\lambda,\zeta})
-\tilde v_{2*}(t^{*\epsilon,\delta,\lambda,\zeta},
y_2^{*\epsilon,\delta,\lambda,\zeta})
\\
&-&(\rho_i-1)y_2^{*\epsilon,\delta,\lambda,\zeta}
\Big(\frac{2}{\epsilon}
(y_1^{*\epsilon,\delta,\lambda,\zeta}
-y_2^{*\epsilon,\delta,\lambda,\zeta})- \delta(\gamma+1)
\exp{( \lambda(T-t^{*\epsilon,\delta,\lambda,\zeta}))}
(y_2^{*\epsilon,\delta,\lambda,\zeta})^\gamma\nonumber\\
&+&
  \zeta(\gamma+1)\frac{1}{(y_2^{*\epsilon,\delta,\lambda,\zeta})^{\gamma
  +2}}
\Big)
\Big)
+y_2^{*\epsilon,\delta,\lambda,\zeta}
\left(
\alpha-\beta+(\beta-\Sum_{i=1}^{d}\rho_i\delta_i\pi_i)_+
\right),\nonumber
\end{eqnarray*}
From the definition of $f$, we have
\begin{eqnarray}\label{t1}
T_2=f(\rho^{*\epsilon,\delta,\lambda,\zeta})&\leq& f({\bf 1})=y_2^{*\epsilon,\delta,\lambda,\zeta}
\big(
\alpha-\beta+(\beta-\Sum_{i=1}^{d}\delta_i\pi_i)_+
\big)\nonumber \\
&\leq & C^{\delta,\zeta},
\end{eqnarray}
where ${\bf 1} $ denotes a $\R^d$-valued vector with all components equal to $1$. Since 
\begin{eqnarray*}
\Phi(t^{*\epsilon,\delta,\lambda,\zeta},
y_1^{*\epsilon,\delta,\lambda,\zeta},
y_2^{*\epsilon,\delta,\lambda,\zeta})
\geq \Phi(t^{*\epsilon,\delta,\lambda,\zeta},
\rho_i^{*\epsilon,\delta,\lambda,\zeta} y_1^{*\epsilon,\delta,\lambda,\zeta},
\rho_i^{*\epsilon,\delta,\lambda,\zeta} y_2^{*\epsilon,\delta,\lambda,\zeta})
\end{eqnarray*}
for all $1\leq i\leq d$, inequality
\reff{uniqueness1} implies
\begin{eqnarray}\label{uniqueness11}
&-&\delta  \lambda
\left(\exp{(\lambda(T-t^{*\epsilon,\delta,\lambda,\zeta}))}
(y_1^{*\epsilon,\delta,\lambda,\zeta})^{\gamma+1}
+\exp{(\lambda(T-t^{*\epsilon,\delta,\lambda,\zeta}))}
(y_2^{*\epsilon,\delta,\lambda,\zeta})^{\gamma+1}\right)\nonumber
\geq T_2-T_1\geq T^{\rho^{*\epsilon,\delta,\lambda,\zeta}},\nonumber
\end{eqnarray}
where
\begin{eqnarray*}
&T^\rho&:=\left\{
\Sum_{i=1}^d 
\pi_i
\Big( (2(\rho_i-1)+1-\rho_i^2)\frac{(
y_1^{*\epsilon,\delta,\lambda,\zeta}-
y_2^{*\epsilon,\delta,\lambda,\zeta})^2}{\epsilon}\right.\nonumber\\
&+&\delta
\big(
(\gamma+1)(\rho_i-1)-\rho_i^{\gamma+1}+1
\big)
\Big(
\exp{( \lambda(T-t^{*\epsilon,\delta,\lambda,\zeta}))}
(y_1^{*\epsilon,\delta,\lambda,\zeta})^{\gamma+1}
\nonumber\\
&+&
\exp{( \lambda(T-t^{*\epsilon,\delta,\lambda,\zeta}))}
(y_2^{*\epsilon,\delta,\lambda,\zeta})^{\gamma+1}
\Big)\nonumber\\
&+&\zeta
\left((\gamma+1)(1-\rho_i)-\frac{1}{\rho_i^{\gamma+1}}+1
\right)
\left(
\frac{1}{(y_1^{*\epsilon,\delta,\lambda,\zeta})^{\gamma+1}}+
\frac{1}{(y_2^{*\epsilon,\delta,\lambda,\zeta})^{\gamma+1}}
\right)
\Big)\nonumber\\
&+&\left.
(y_2^{*\epsilon,\delta,\lambda,\zeta}-y_1^{*\epsilon,\delta,\lambda,\zeta})
\big(
\alpha-\beta+(\beta-\Sum_{i=1}^{d}\rho_i\delta_i\pi_i)_+
\big)
\Big)
\right\}.
\end{eqnarray*}
Sending $\epsilon \longrightarrow 0^+$ ,
we obtain 
\begin{eqnarray}\label{uniqueness111}
&-&2 \delta \lambda \exp{\left( \lambda(T-t^{*\delta,\lambda,\zeta})\right)}
(
y^{*\delta,\lambda,\zeta})^{\gamma+1}\nonumber\\
&\geq &
\Sum_{i=1}^d 
\pi_i
\Big(
2 \delta
\left((\gamma+1)(\rho_i^{*\delta,\lambda,\zeta}-1)+
( \rho_i^{*\delta,\lambda,\zeta})^{\gamma+1}-1
\right)
\exp{\left(\lambda(T-t^{*\delta,\lambda,\zeta})\right)}
(y^{*\delta,\lambda, \zeta})^{\gamma+1}\nonumber\\
&+&2 \zeta\left((\gamma+1)(1- \rho_i^{*\delta,\lambda,\zeta})-
\frac{1}{(\rho_i^{*\delta,\lambda,\zeta})^{\gamma+1}}+1\right)
\frac{1}{(y^{*\delta,\lambda,\zeta})}\Big).
\end{eqnarray}
Using the boundedness of $y^{*\delta,\lambda,\zeta}$ and $\rho^{*\delta,\lambda,\zeta}\in \Sigma^{'}$, 
along a subsequence $(y^{*\delta,\lambda,\zeta},\rho^{*\delta,\lambda,\zeta}) $ 
converges when $\lambda \longrightarrow \infty$. Let's denote 
$(y^{*\delta,\zeta},\rho^{*\delta,\zeta} )$ its limit.
Sending $\lambda \longrightarrow \infty$ in
inequality \reff{uniqueness111}, we obtain 
\begin{eqnarray*}
& &\Sum_{i=1}^d 
\pi_i
\Big(
2 \delta
\left((\gamma+1)(\rho_i^{*\delta,\zeta}-1)+
( \rho_i^{*\delta,\zeta})^{\gamma+1}-1
\right)
\exp{\left(\lambda(T- t^{*\delta,\zeta})\right)}
(y^{*\delta, \zeta})^{\gamma+1}\nonumber\\
&+&2 \zeta\left((\gamma+1)(1- \rho_i^{*\delta,\zeta})-
\frac{1}{( \rho_i^{*\delta,\zeta})^{\gamma+1}}+1\right)
\frac{1}{y^{*\delta,\zeta}}\Big)=-\infty
\end{eqnarray*}
which implies, there exists $i_0$, $1\leq i_0\leq d$ such that $\rho_{i_0}^{*\delta,\zeta}=0$. Sending $\epsilon \longrightarrow 0$ and $\lambda \longrightarrow \infty$ in inequality \reff{t1}, we obtain $f(\rho^{*\delta,\zeta})=\infty\leq C^{\delta,\zeta}$ which is false. \\
Sending $\epsilon\longrightarrow \infty$ in inequality \reff{uniqueness2}, we
have
\begin{eqnarray*}
-\delta \exp{\left( \lambda (T-
 t^{*\delta,\lambda,\zeta})\right)}
(y^{*\delta,\lambda,\zeta})+
\zeta \frac{1}{(y^{*\delta, \lambda,\zeta})^{\gamma+2}}\geq 0,
\end{eqnarray*}
which implies
\begin{eqnarray}\label{2terme}
\frac{\delta \exp{\left(\lambda(T-
 t^{*\delta,\lambda,\zeta})\right)}}{\zeta}
\leq \frac{1}{(M_1^{\delta,\zeta})^{2(\gamma+1)}}.
\end{eqnarray}
Using the boundedness of $t^{*\delta,\lambda,\zeta}$,
along a subsequence $ t^{*\delta,\lambda,\zeta}$ 
converges when $\lambda \longrightarrow \infty$. From inequality \reff{2terme}, we have necessarily 
$t^{*\delta,\zeta}=T$.\\
From inequality
\begin{eqnarray*}
\Phi(t,y,y)\leq
\Phi(t^{*\epsilon,\delta,\lambda,\zeta},
y_1^{*\epsilon,\delta,\lambda,\zeta},
y_2^{*\epsilon,\delta,\lambda,\zeta})
\end{eqnarray*} 
and since $\Phi$ is upper semi-continuous, we deduce
that 
\begin{eqnarray*}
\Phi(t,y,y)&\leq& 
\Limsup_{ \lambda \rightarrow \infty} \Limsup_{ \epsilon \rightarrow 0}
\Phi(t^{*\epsilon,\delta,\lambda,\zeta},
y_1^{*\epsilon,\delta,\lambda,\zeta},
y_2^{*\epsilon,\delta,\lambda,\zeta})\\
&\leq&\Phi(T,
y^{\delta,\zeta},
y^{\delta,\zeta})\\
&\leq & \tilde v^*_1(T,y^{\delta,\zeta})
-\tilde v_{2*}(T,y^{\delta,\zeta})\leq 0.
\end{eqnarray*}
Sending 
$\delta$, $\zeta$ $\longrightarrow 0^+$, we obtain
\begin{eqnarray*}
\tilde v^*_1(t,y) \leq \tilde v_{2*}(t,y),\,\mbox{ for all }(t,y)\in [0,T]\times (0,\infty).
\end{eqnarray*} 
 and so Theorem \ref{comparaison} is proved. 
\ep

\end{document}